\documentclass[12pt]{amsart}
\usepackage[vmargin=2.6cm,hmargin=2.2cm]{geometry} 
\usepackage[colorlinks,citecolor=blue,urlcolor=blue,bookmarks=false,hypertexnames=true]{hyperref}
\usepackage{color}
\usepackage{amsthm}
\usepackage{graphicx}
\usepackage{times}
\usepackage{xcolor}
\graphicspath{ {./images/} }
\usepackage{tikz}
\usepackage{multirow}
\usepackage{caption}
\usepackage{subcaption}
\usepackage[square,numbers]{natbib}
\usepackage{comment}
\usepackage[title]{appendix}
\usepackage{fancyhdr}	

\allowdisplaybreaks[3] 

\numberwithin{equation}{section}
\newtheorem{theorem}{Theorem}[section]

\newtheorem{lemma}{Lemma}[section]

\newtheorem{remark}{Remark}[section]

\newtheorem{assumption}{Assumption}[section]
\newtheorem{definition}{Definition}[section]

\newcommand{\dd}{\operatorname{d}\! }
\newcommand{\dt}{\operatorname{d}\! t}
\newcommand{\de}{\operatorname{d}\! e}
\newcommand{\ds}{\operatorname{d}\! s}

\newcommand{\dw}{\operatorname{d}\! W}

\newcommand{\cM}{\ensuremath{\mathcal{M}}}

\newcommand{\cE}{\ensuremath{\mathcal{E}}}
\newcommand{\argmin}{\ensuremath{\operatorname*{argmin}}}

\newcommand{\esssup}{\ensuremath{\operatorname*{ess\;sup}}}

\newcommand{\dualpi}{\widehat{\Pi}}
\newcommand{\ep}{\varepsilon}
\newcommand{\lnu}{\ensuremath{L^{2,\nu}}}

\newcommand{\ltwonu}{\ensuremath{L^{2, \nu}_{\mathcal{P}}}}
\newcommand{\linnu}{\ensuremath{L^{\infty,\nu}_{\mathcal{P}}}}
\newcommand{\nn}{\nonumber}

\newcommand{\E}{{\mathbb E}}
\newcommand{\R}{{\mathbb R}}
\newcommand{\pf}{\noindent\textbf{Proof:} }
\newcommand{\eof}{\hfill{$\Box$}}

\newcommand{\dpt}{\dd\mathbb{P}\otimes \dt\textrm{-a.e.}}

\title[Comparison theorems for multi-dimensional BSDEs with jumps]{Comparison theorems for multi-dimensional BSDEs with jumps and applications to constrained stochastic linear-quadratic control}
\author[Hu]{Ying Hu}
\author[Shi]{Xiaomin Shi}
\author[Xu]{Zuo Quan Xu}
\date{\today}

\address{Y.~Hu: Univ Rennes, CNRS, IRMAR-UMR 6625, F-35000 Rennes, France.}
\email{{ying.hu@univ-rennes1.fr}}
\address{X.~Shi: School of Statistics and Mathematics, Shandong University of Finance and Economics, Jinan
250100, China.}
\email{{shixm@mail.sdu.edu.cn}}
\address{Z.~Xu: Department of Applied Mathematics, The Hong Kong Polytechnic University, Kowloon, Hong Kong.}
\email{{maxu@polyu.edu.hk}}

\begin{document}
\setcitestyle{numbers}
\maketitle


\begin{abstract}
In this paper, we, for the first time, establish two comparison theorems for multi-dimensional backward stochastic differential equations with jumps. Our approach is novel and completely different from the existing results for one-dimensional case.
Using these and other delicate tools, we then construct solutions to coupled two-dimensional stochastic Riccati equation with jumps in both standard and singular cases.
In the end, these results are applied to solve a cone-constrained stochastic linear-quadratic and a mean-variance portfolio selection problem with jumps. Different from no jump problems, the optimal (relative) state processes may change their signs, which is of course due to the presence of jumps.
\end{abstract}

\subsection*{Keywords:} Backward stochastic differential equations with jumps, multi-dimensional comparison theorem, stochastic Riccati equation with jumps, cone-constrained linear-quadratic control, mean-variance problem.
\subsection*{Mathematics Subject Classification (2020):} 60H30. 60J76. 93E20. 91G10.

\section{Introduction}
The study of backward stochastic differential equations (BSDEs, for short) can be dated back to Bismut \cite{Bismut}, who studied the linear case, as an adjoint equation in the Pontryagin stochastic maximum principle. The general Lipschitz continuous case was later resolved in the seminal paper of Pardoux and Peng \cite{PP}. Since then, BSDEs have attracted strong interest of many researchers and found widely applications in partial differential equations, stochastic control, stochastic differential game and mathematical finance; see, e.g., \cite{DP,EH,EPQ, EPQ01, HP, HZ, Peng}.
In particular, the solvability of quadratic BSDEs in one-dimensional case was firstly obtained in Kobylanski \cite{Ko}, and then generalized to multi-dimensional case by \cite{FHT,HT,Luo,Te}.

BSDEs that are driven by a Brownian motion and a Poisson random measure, which are named as BSDE with jumps (BSDEJ) in this paper, were firstly tackled by Tang and Li \cite{TL}, then followed notably by Barles, Buckdahn and Pardoux \cite{BBP}, Royer \cite{Royer}, Quenes and Sulem \cite{QS} in the Lipschitz case. Quadratic BSDEs with jumps and their applications in utility maximization problems have also been investigated; see, e.g., Antonelli and Mancini \cite{AM}, Kazi-Tani, Possama\"{\i} and Zhou \cite{KPZ}, Laeven and Stadje \cite{LS}, Morlais \cite{Mo, Mo2} among many others. Please refer to Papapantoleon, Possama\"{\i}, Saplaouras \cite{PPS} for a synopsis of these topics.

BSDEs arising from stochastic linear quadratic (LQ) control problems, called the stochastic Riccati equations (SREs), form an important class of BSDEs.
In these BSDEs, the first unknown variable appears on the denominator
and the second unknown variable grows quadratically in the generator. These features distinguish them from those well-studied BSDEs with Lipschitz or quadratic growth generators, so that they have to be studied by new methods.

Bismut \cite{Bismut76} firstly found that a linear state feedback form optimal control for a stochastic LQ control problem is available, provided that its associated SRE admits a solution in some suitable space. Unfortunately he could not show the existence of such a solution in general.
Nowadays numerous progresses have been made in solving SREs. Kohlmann and Tang \cite{KT} resolved the existence and uniqueness issues for one-dimensional SREs, then Tang \cite{Tang03,Tang15} resolved the matrix-valued case using the stochastic maximum principle and dynamic programming method respectively. Sun, Xiong and Yong \cite{SXY} studied the indefinite case. Inspired by Tang's \cite{Tang15} dynamic programming method, Zhang, Dong and Meng \cite{ZDM} established the existence and uniqueness of solutions to SREs with jumps (SREJ). Li, Wu and Yu \cite{LWY} studied the indefinite case using a so-called relax compensator.

Motivated by the mean-variance (MV) portfolio selection problem with no-shorting constraints, Hu and Zhou \cite{HZ} studied cone-constrained stochastic LQ problem and found that the optimal control takes a piecewise ($0$ is the unique segment point) linear state feedback form. The associated SRE is a two-dimensional, but decoupled, BSDE. Hence it can be treated separately as two one-dimensional BSDEs. The solvability was established with the aid of quadratic BSDE theory and truncation techniques. The decoupling phenomenon lies in the fact that the optimal state process will not change its sign (namely not cross $0$), i.e. it will stay positive (resp. negative) if the initial state is positive (resp. negative). Dong \cite{Dong} generalized the model in \cite{HZ} to incorporate a jump by the enlargement of filtration framework. The corresponding SRE is a coupled two-dimensional BSDEJ, whose solvability is obtained by solving two recursive systems of BSDEs driven only by Brownian motions. This decomposition approach works only in the filtration enlargement theory; see also Kharroubi, Lim and Ngoupeyou \cite{KLN}, Hu, Shi and Xu \cite{HSX2} for the unconstrained or regime switching case.
Czichowsky and Schweizer \cite{CS} extended the cone-constrained MV model to a general semi-martingale framework, but they can not solve the two-dimensional SREJ. They claimed that ``finding a solution by general BSDE techniques seems a formidable challenge" in \cite[Remark 4.8]{CS}.

This paper is intended as an attempt to cope with the formidable challenge indicated in \cite{CS}. Our main contribution is to resolve the solvability of a two-dimensional coupled SREJ in the Wiener-Poisson world via pure BSDE techniques. Although one can consider the more general semi-martingale framework, we will focus on the Wiener-Poisson world as SREJ in this case takes more concrete structures for presentation and illustration. We establish the solvability for both standard and singular cases, containing the SREJ emerging in the cone-constrained MV problem as an special example. Since the existing approximation procedures in Kohlmann and Tang \cite{KT} and our previous work \cite{HSX} can not be applied to the present problem, we provide a new approximation procedure to achieve the goal.

A crucial and novel tool used in the approximation approach is our new comparison theorems for BSDEJs. We establish two comparison theorems which seem to be the first ones for multi-dimensional case. The first one requires a locally Lipschitz condition for one generator
(see Remark \ref{remark:cond-growth}) and works for bounded state processes, whereas the second one requires the globally Lipschitz condition for both generators and works for square integrable state processes.

Most of existing comparison theorems for BSDEJs require the condition $\gamma>-1$ (see Remark \ref{remark:cond-gamma}) or even stronger $\gamma>-1+\ep$ in order to utilize the Girsanov theorem; see, e.g., Barles, Buckdahn and Pardoux \cite{BBP} and Royer \cite{Royer}. To the best of our knowledge, Quenez and Sulem's \cite{QS} comparison theorem is the only one that relaxes the condition to $\gamma\geq-1$. Without resort to the Girsanov theorem, they used the conditional expectation representation of one-dimensional linear BSDEJs to establish their comparison theorem. Nevertheless all of these existing comparison theorems for BSDEJs can only deal with one-dimensional case. In our approximation procedure, however, the SREJ is a fully coupled two-dimensional BSDEJ, therefore comparison theorems for multi-dimensional BSDEJs are strongly appealing. It is worth pointing that the conditional expectation representation method used in \cite{QS} cannot be applied to multi-dimensional BSDEJs.
In this paper, we propose a completely different approach to establish our comparison theorems for multi-dimensional BSDEJs for the first time. We achieve the goal by directly analyzing $((\delta Y_t)^+)^2$ with the aid of the Meyer-It\^o formula and utilizing a tricky elementary inequality (Lemma \ref{Ausefullinequality}) that works for $\gamma\geq-1$.
Note one cannot expect to extend the results to the case $\gamma<-1$ since counter-examples do exist in this case; see \cite[Remark 2.7]{BBP}.

With the help of the new comparison theorems for multi-dimensional BSDEJs, we can construct solutions to the two-dimensional coupled SREJ in both standard and singular cases. We then apply the result to solve a cone-constrained stochastic LQ problem with jumps and obtain the efficient portfolio for a MV problem with jumps. It is worth pointing out that even without the cone-constraint, MV problems with jumps have not been investigated thoroughly. Lim \cite{Lim} studied such a problem, but he assumed all the coefficients are predictable with respect to the Brownian filtration, rendered the corresponding SRE exactly the same as that in the model without jumps. On the other hand, Zhang, Dong and Meng \cite{ZDM} examined stochastic LQ problems with jumps, but they assumed the control weight in the running cost is uniformly positive so that their result cannot solve the corresponding MV problem where the control weight is 0. We will not only solve the MV problem with jumps, but also incorporate convex cone-constraint, especially covering the famous no-shorting constraints. By adding cone-constraint, the associated SREJ becomes a fully coupled two-dimensional BSDEJ, thus causing notably nontrivial difficulty in its solvability.

The rest part of this paper is organized as follows. Section \ref{section:comp} is devoted to proving two comparison theorems for multi-dimensional BSDEJs.
In Section \ref{section:LQ}, we study a cone-constrained stochastic LQ control problem with jumps and prove the existence and uniqueness of solution to the associated SREJ.
In Section \ref{section:MV}, we solve a cone-constrained MV problem. Appendix \ref{appnA} provides a heuristical derivation of the SREJ. A lengthy and complementary proof of Theorem \ref{Th:standard} is relegated to Appendix \ref{appnB}.

\section{Comparison theorems for multi-dimensional BSDEJs}\label{section:comp}
Let $(\Omega, \mathcal F, \mathbb{F}, \mathbb{P})$ be a fixed complete filtered probability space. The filtration $\mathbb{F}=\{\mathcal F_t, t\geq0\}$ is generated by two independent random sources augmented by all $\mathbb{P}$-null sets: one is
a standard $n$-dimensional Brownian motion $W_t=(W_{1,t}, \ldots, W_{n,t})^{\top}$, and the other one is
a Poisson random measure $N(\dt,\de)$ defined on $\R_+\times\cE$ induced by a stationary Poisson point process with a stationary compensator (intensity measure) given by $\nu(\de)\dt$ satisfying $\nu(\cE)<\infty$, where $\mathcal{\cE}\subseteq \R^{\ell}\setminus\{0\}$ is a nonempty Borel subset of the $\ell$-dimensional Euclidean space $\R^{\ell}$.
We use an increasing sequence $\{T_n\}_{n\in\mathbb{N}}$ to denote the jump times of underlying Poisson point process.
The compensated Poisson random measure is denoted by $ \widetilde N(\dt,\de)$.
For the ease of notations, we only consider one-dimensional Poisson random measure, although the results of this paper can be generalized to the multi-dimensional case without essential difficulties. Throughout the paper, let $T$ denote a fixed positive constant, $\mathcal{P}$ denote the $\mathbb{F}$-predictable $\sigma$-field on $\Omega\times[0,T]$, and $\mathcal{B}(\cE)$ denote the Borel $\sigma$-algebra of $\cE$.

We denote by $\R^\ell$ the set of $\ell$-dimensional column vectors, by $\R^\ell_+$ the set of vectors in $\R^\ell$ whose components are nonnegative, by $\R^{\ell\times n}$ the set of $\ell\times n$ real matrices, and by $\mathbb{S}^n$ the set of symmetric $n\times n$ real matrices. Therefore, $\R^\ell\equiv\R^{\ell\times 1}$. For any vector $Y$, we denote $Y_i$ as its $i$-th component.
For any matrix $M=(m_{ij})$, we denote its transpose by $M^{\top}$, and its norm by $|M|=\sqrt{\sum_{ij}m_{ij}^2}$. If $M\in\mathbb{S}^n$ is positive definite (resp. positive semidefinite), we write $M>$ (resp. $\geq$) $0.$ We write $A>$ (resp. $\geq$) $B$ if $A, B\in\mathbb{S}^n$ and $A-B>$ (resp. $\geq$) $0.$ We use the standard notations $x^+=\max\{x, 0\}$ and $x^-=\max\{-x, 0\}$ for $x\in\R$ and define a set $\cM=\{1,2,...,\ell\}$.
We will use the elementary inequality $|a^{\top}b|\leq c|a|^2+\frac{|b|^2}{2c}$ for any $a,b\in\R^{n}, c>0$ frequently throughout the paper without claim.

We use the following spaces throughout the paper:
\vspace{-5pt}
\begin{align*}
L^2_{\mathcal{F}_T}(\Omega;\R)&=\Big\{\xi:\Omega\to \R
\;\Big|\;\mbox{$\xi$ is $\mathcal{F}_{T}$-measurable, and $\E|\xi|^{2}<\infty$}\Big\}, \\
L^{\infty}_{\mathcal{F}_T}(\Omega;\R)&=\Big\{\xi:\Omega\to \R
\;\Big|\;\mbox{$\xi$ is $\mathcal{F}_{T}$-measurable, and essentially bounded}\Big\}, \\
L^{2}_{\mathbb F}(0, T;\R)&=\Big\{\phi:\Omega\times[0, T]\to \R
\;\Big|\;\mbox{$\phi$ is $\mathcal{P}$-measurable and $\E\int_{0}^{T}|\phi_t|^{2}\dt<\infty$}\Big\}, \\
L^{\infty}_{\mathbb{F}}(0, T;\R)&=\Big\{\phi:\Omega\times[0, T]\to \R
\;\Big|\;\mbox{$\phi$ is $\mathcal{P}$-measurable and essentially bounded} \Big\},\\
\lnu(\R)&=\Big\{\phi:\cE\to \R
\;\Big|\;\mbox{$\phi$ is $\mathcal{B}(\cE)$-measurable and $||\phi||^2_{\nu}:=\int_{\cE} |\phi(e)|^2\nu(\de)<\infty$}\Big\},\\
L^{\infty,\nu}(\R)&=\Big\{\phi:\cE\to \R
\;\Big|\;\mbox{$\phi$ is $\mathcal{B}(\cE)$-measurable and essentially bounded w.r.t. $\dd\nu$}\Big\}, \\
\ltwonu(0, T;\R)&=\Big\{\phi:\Omega\times[0, T]\times \cE\to \R
\;\Big|\;\mbox{$\phi$ is $\mathcal{P}\otimes\mathcal{B}(\cE)$-measurable }\\
&\qquad\quad\ \mbox{and $\E\int_{0}^{T}\int_{\cE}|\phi_t(e)|^{2} \nu(\de)\dt<\infty$}
\Big\},\\
\linnu(0, T;\R) &=\Big\{\phi:\Omega\times[0, T]\times \cE\to \R
\;\Big|\;\phi\mbox{ is $\mathcal{P}\otimes\mathcal{B}(\cE)$-measurable and}\\
&\qquad\quad\mbox{ essentially bounded w.r.t. $\dd\mathbb{P}\otimes \dt\otimes \dd\nu$}\Big\},\\
S^{2}_{\mathbb{F}}(0,T;\R)&=\Big\{\phi:\Omega\times[0,T]\to \R
\;\Big|\;(\phi_{t})_{0\leq t\leq T}\mbox{ is c\`ad-l\`ag, $\mathbb{F}$-adapted}\\
&\qquad\quad \mbox{ and $\E\sup_{0\leq t\leq T}|\phi_t|^2<\infty$}\Big\},\\
S^{\infty}_{\mathbb{F}}(0,T;\R)&=\Big\{\phi:\Omega\times[0,T]\to \R
\;\Big|\;(\phi_{t})_{0\leq t\leq T}\mbox{ is c\`ad-l\`ag, $\mathbb{F}$-adapted}\\
&\qquad\quad \mbox{ and essentially bounded}\Big\}.
\end{align*}
These definitions are generalized in the obvious way to the cases that $\R$ is replaced by $\R^n$, $\R^{n\times \ell}$ or $\mathbb{S}^n$.
Arguments $s$, $t$ and $\omega$, or statements ``almost surely'' (a.s.) and ``almost everywhere'' (a.e.), may be suppressed for simplicity in many circumstances when no confusion occurs. We shall use $c$ to represent a generic positive constant which can be different from line to line. All the equations and inequalities in subsequent analysis shall be understood in the sense that $\dd\mathbb{P}$-a.s. or $\dd\nu$-a.e. or $\dpt$ or $\dd\mathbb{P}\otimes \dt\otimes \dd\nu$-a.e. etc.

In this paper, any $\ell$-dimensional backward stochastic differential equation with jumps (BSDEJ) (on $[0,T]$) is characterized by a pair $(\xi, f)$, in which $\xi:\Omega\to \R^{\ell}$ is called the terminal value which is an $\mathcal{F}_{T}$-measurable random vector, and $f:\Omega\times[0,T]\times\R^\ell\times\R^{n\times \ell}\times L^{2,\nu}(\R^\ell)\to \R^\ell$ is called the generator which is a $\mathcal{P}\otimes\mathcal{B}(\R^\ell)\otimes\mathcal{B}(\R^{n\times\ell})\otimes\mathcal{B}(L^{2,\nu}(\R^\ell))$-measurable process.
We call the BSDEJ $\ell$-dimensional as its state process is $\R^\ell$-valued. We often rewrite in its component form for the ease of presentations.

%
%

\subsection{Comparison theorem for bounded processes}

We first prove a comparison theorem where the state processes are essentially bounded.

\begin{theorem}
\label{comparisonbound}
Suppose, for every $ i\in \cM$,
\begin{align*}
(Y_i, Z_i, \Phi_i), (\overline Y_i, \overline Z_i, \overline\Phi_i)\in S^{\infty}_{\mathbb{F}}(0,T;\R)\times L^{2}_{\mathbb F}(0, T;\R^{n})\times L^{2,\nu}_{\mathcal{P}}(0, T;\R),
\end{align*}
and they satisfy BSDEJs
\begin{align}\label{CY}
Y_{i,t}&=\xi_i+\int_t^T f_i(s,Y_{s-}, Z_{i,s}, \Phi_{s})\ds\nn\\
&\qquad-\int_t^T Z_{i,s}^{\top}\dw_s-\int_t^T\int_{\cE}\Phi_{i,s}(e) \widetilde N(\ds,\de)~\dpt,
\end{align}
and
\begin{align}\label{CYoverline}
\overline Y_{i,t}&=\overline\xi_i+\int_t^T \overline f_i(s, \overline Y_{s-}, \overline Z_{i,s}, \overline \Phi_{s})\ds\nn\\
&\qquad-\int_t^T\overline Z_{i,s}^{\top}\dw_s-\int_t^T\int_{\cE}\overline\Phi_{i,s}(e) \widetilde N(\ds,\de)~\dpt.
\end{align}
Also suppose that, for all $i\in\cM$ and $s\in[0,T]$,
\begin{enumerate}

\item \label{cond-boundary}
$\xi_i\leq\overline\xi_i$;

\item \label{cond-gamma}
there exists a constant $c>0$ such that
\begin{align*}
&\quad\;f_i(s,Y_{s-}, Z_{i,s}, \Phi_{1,s}, \cdots, \Phi_{i,s}, \cdots, \Phi_{\ell,s})\\
&\qquad\quad-f_i(s,Y_{s-}, Z_{i,s}, \Phi_{1,s}, \cdots, \overline\Phi_{i,s}, \cdots, \Phi_{\ell,s}) \\
&\leq c \int_{\cE} (\Phi_{i,s}(e)-\overline \Phi_{i,s}(e))^{+}\nu(\de)
+\int_{\cE} |\Phi_{i,s}(e)-\overline \Phi_{i,s}(e)|\nu(\de);
\end{align*}


\item \label{cond-growth} there exists a constant $c>0$ such that
\begin{align*}
&\quad\;f_i(s,Y_{s-}, Z_{i,s}, \Phi_{1,s}, \cdots, \overline\Phi_{i,s}, \cdots, \Phi_{\ell,s})-f_i(s, \overline Y_{s-}, \overline Z_{i,s}, \overline \Phi_{s})\\
&\leq c\Big (| Y_{i,s-}- \overline Y_{i,s-}|+\sum_{j\neq i} (Y_{j,s-}- \overline Y_{j,s-})^{+}+|Z_{i,s}- \overline Z_{i,s}|\\
&\qquad\quad+\sum_{j\neq i} \int_{\cE}( Y_{j,s-}+\Phi_{j,s}(e)-\overline Y_{j,s-}
-\overline\Phi_{j,s}(e))^+\nu(\de)\Big); ~\text{and}
\end{align*}

\item \label{cond-size} $f_i(s, \overline Y_{s-}, \overline Z_{i,s}, \overline \Phi_{s})
\leq \overline f_i(s, \overline Y_{s-}, \overline Z_{i,s}, \overline \Phi_{s}).$
\end{enumerate}
Then $Y_{i}\leq \overline Y_{i}$ for all $i\in\cM$.
\end{theorem}
To prove this theorem, we need the following critical elementary result.
\begin{lemma}\label{Ausefullinequality}
For all $(x,y)\in\R\times\R$ and $c\geq -1$, we have
\begin{align*}
[(x+y)^+]^2-(x^+)^2-2(1+c)x^+y\geq -(c^2\vee 1)(x^{+})^2.
\end{align*}
\end{lemma}
\pf
There are three cases:
\begin{itemize}
\item If $x\leq 0$, then
\begin{align*}
[(x+y)^+]^2-(x^+)^2-2(1+c)x^+y=[(x+y)^+]^2\geq 0= -(c^2\vee 1)(x^+)^2.
\end{align*}
\item If $y\leq 0$, then, since $c\geq -1$,
\begin{align*}
[(x+y)^+]^2-(x^+)^2-2(1+c)x^+y\geq [(x+y)^+]^2-(x^+)^2\geq -(x^+)^2\geq -(c^2\vee 1)(x^+)^2.
\end{align*}
\item If $x\geq 0$ and $y\geq0$, then
\begin{align*}
[(x+y)^+]^2-(x^+)^2-2(1+c)x^+y=y^2-2c xy&=(y-c x)^2-c^2x^2\\
&\geq -(c^2\vee 1)x^2=-(c^2\vee 1)(x^{+})^2.\qquad
\end{align*}
\end{itemize}
The proof is complete.
\eof
\bigskip\\

\noindent\textbf{Proof of Theorem \ref{comparisonbound}.}\\
For $t\in[0, T]$ and $i\in\cM$, set
\[\delta Y_{i,t}=Y_{i,t}-\overline Y_{i,t}, \ \delta Z_{i,t}=Z_{i,t}-\overline Z_{i,t}, \ \delta \Phi_{i,t}=\Phi_{i,t}-\overline \Phi_{i,t}.\]
Applying the Meyer-It\^o formula \cite[Chapter IV, Theorem 70]{Protter} to $(\delta Y_{i,t})^+$, we get,
\begin{align*}
\dd\; (\delta Y_{i,t})^+&=-\mathbf{1}_{\{\delta Y_{i,t-}>0\}}[f_i(t,Y_{t-}, Z_{i,t}, \Phi_{t} )-\overline f_i(t, \overline Y_{t-}, \overline Z_{i,t}, \overline \Phi_{t})]\dt\\
&\quad\;+\int_\cE[(\delta Y_{i,t-}+\delta \Phi_{i,t}(e))^+ -(\delta Y_{i,t-}^i)^{+}-\mathbf{1}_{\{\delta Y_{i,t-}>0\}}\delta \Phi_{i,t}(e)]\nu(\de)\dt+\frac{1}{2}\dd L_{i,t}\\
&\quad\;+\mathbf{1}_{\{\delta Y_{i,t-}>0\}}\delta Z_{i,t}^{\top}\dw_t+\int_\cE [(\delta Y_{i,t-}+\delta \Phi_{i,t}(e))^+ -(\delta Y_{i,t-})^+] \widetilde N(\dt,\de),
\end{align*}
where $L_{i,t}$ is the local time of $\delta Y_{i,t}$ at $0$.
Since $\delta Y_{i,t-} \dd L_{i,t}=0$, applying It\^{o}'s formula to $((\delta Y_{i,t})^+)^2$ yields
\begin{align}\label{deltaYsquare}
\dd\; ((\delta Y_{i,t})^+)^2
&=-2(\delta Y_{i,t-})^+[f_i(t,Y_{t-}, Z_{i,t}, \Phi_{t} )-\overline f_i(t, \overline Y_{t-}, \overline Z_{i,t}, \overline \Phi_{t})]\dt+\mathbf{1}_{\{\delta Y_{i,t-}>0\}} |\delta Z_{i,t}|^2\dt\nn\\
&\quad\;+\int_\cE [((\delta Y_{i,t-}+\delta \Phi_{i,t}(e))^+)^2 -((\delta Y_{i,t-})^+)^2-2(\delta Y_{i,t-})^+\delta \Phi_{i,t}(e)]\nu(\de)\dt\nn\\
&\quad\;+2(\delta Y_{i,t-})^+\delta Z_{i,t}^{\top}\dw_t+\int_{\cE}[((\delta Y_{i,t-}+\delta \Phi_{i,t}(e))^+)^2 -((\delta Y_{i,t-})^+)^2] \widetilde N(\dt,\de).
\end{align}
Using the condition \ref{cond-size} and inserting two zero-sum terms, we get
\begin{align*}
&\quad\; f_i(t,Y_{t-}, Z_{i,t}, \Phi_{t} )-\overline f_i(t, \overline Y_{t-}, \overline Z_{i,t}, \overline \Phi_{t})\nn\\
&\leq f_i(t,Y_{t-}, Z_{i,t}, \Phi_{t} )- f_i(t, \overline Y_{t-}, \overline Z_{i,t}, \overline \Phi_{t})\nn\\
&=[f_i(t,Y_{t-}, Z_{i,t}, \Phi_{t} ) -f_i(t,Y_{t-}, Z_{i,t}, \Phi_{1,t}, \cdots, \overline\Phi_{i,t}, \cdots, \Phi_{\ell,t})]\nn\\
&\quad\;+[f_i(t,Y_{t-}, Z_{i,t}, \Phi_{1,t}, \cdots, \overline\Phi_{i,t}, \cdots, \Phi_{\ell,t})-f_i(t, \overline Y_{t-}, \overline Z_{i,t}, \overline \Phi_{t})].
\end{align*}
By the conditions \ref{cond-gamma}, the first difference on the right hand side (RHS) in above is upper bounded by
\[ \int_{\cE}\gamma_{i,t}(e) \delta\Phi_{i,t}(e)\nu(\de),\]
where
\begin{align*}
\gamma_{i,t}(e)=\begin{cases}
c&\quad \text{ if }\delta\Phi_{i,t}(e)\geq 0;\\
-1&\quad \text{ if }\delta\Phi_{i,t}(e)<0.
\end{cases}
\end{align*}
By the conditions \ref{cond-growth}, the second difference on the RHS is upper bounded respectively by
\[c\Big (|\delta Y_{i,t-}|+\sum_{j\neq i} (\delta Y_{j,t-})^{+}+|\delta Z_{i,t}|+\sum_{j\neq i} \int_{\cE}(\delta Y_{j,t-}+\delta\Phi_{j,t}(e))^+\nu(\de)\Big).\]
%
Using these estimates and $\nu(\cE)<\infty$, we deduce that
\begin{align}\label{deltaf}
&2(\delta Y_{i,t-})^+[f_i(t,Y_{t-}, Z_{i,t}, \Phi_{t} )-\overline f_i(t, \overline Y_{t-}, \overline Z_{i,t}, \overline \Phi_{t})] \nn\\
&\leq c\sum_{i=1}^{\ell}((\delta Y_{i,t-})^+)^2+\mathbf{1}_{\{\delta Y_{i,t-}>0\}}|\delta Z_{i,t}|^2+\mathbf{1}_{\{\delta Y_{i,t-}>0\}}\sum_{j\neq i} \int_{\cE}((\delta Y_{j,t-}+\delta\Phi_{j,t}(e))^+)^{2}\nu(\de)\nn\\
&\quad\;+2(\delta Y_{i,t-})^+\int_{\cE}\gamma_{i,t}(e) \delta\Phi_{i,t}(e)\nu(\de).
\end{align}
Integrating from $t$ to $T$ in \eqref{deltaYsquare}, taking conditional expectation and using \eqref{deltaf}, we obtain
\begin{align*}
&\quad((\delta Y_{i,t})^+)^2\\
&\leq \mathbb{E}_t\int_t^T \Big(c\sum_{i=1}^{\ell}((\delta Y_{i,s-})^+)^2+\mathbf{1}_{\{\delta Y_{i,s-}>0\}}\sum_{j\neq i}\int_{\cE}((\delta Y_{j,s-}+\delta \Phi_{j,s}(e))^+)^2\nu(\de)\Big)\ds\\
&-\mathbb{E}_t\int_t^T\int_\cE \Big[((\delta Y_{i,s-}+\delta \Phi_{i,s}(e))^+)^2 -((\delta Y_{i,s-})^+)^2-2(1+\gamma_{i,s}(e))(\delta Y_{i,s-})^+\delta \Phi_{i,s}(e)\Big]\nu(\de)\ds.
\end{align*}
Because $\gamma_{i}\in \linnu(0, T;\R)$ and $\gamma_{i}\geq -1$, it follows from Lemma \ref{Ausefullinequality} that
\begin{multline}\label{deltaPhi}
-[((\delta Y_{i,s-}+\delta \Phi_{i,s}(e))^+)^2 -((\delta Y_{i,s-})^+)^2-2(1+\gamma_{i,s}(e))(\delta Y_{i,s-})^+\delta \Phi_{i,s}(e)] \\
\qquad\leq (\gamma_{i,s}(e)^{2}\vee 1) ((\delta Y_{i,s-})^+)^2 \leq c ((\delta Y_{i,s-})^+)^2.
\end{multline}
Combining the above estimates and using $\nu(\cE)<\infty$, we obtain
\begin{align}\label{deltaY}
((\delta Y_{i,t})^+)^2
&\leq c\mathbb{E}_t\int_t^T \sum_{i=1}^{\ell}((\delta Y_{i,s-})^+)^2\ds +\sum_{j\neq i}\E_t\int_t^T\int_{\cE}((\delta Y_{j,s-}+\delta \Phi_{j,s}(e))^+)^2\nu(\de)\ds,
\end{align}
where the constant $c$ is independent of $t$, $T$ and $i$.

Note that
\begin{align}\label{jumpterm}
\E_t\int_t^T\int_{\cE}((\delta Y_{j,s-}+\delta \Phi_{j,s})^+)^2\nu(\de)\ds&=\E_t\int_t^T\int_{\cE}((\delta Y_{j,s-}+\delta \Phi_{j,s}(e))^+)^2N(\ds,\de)\nn\\
&=\E_t\Big[\sum_{n\in\mathbb{N}, \ t<T_n\leq T}((\delta Y_{j,T_n-}+\delta \Phi_{j,T_n}(\Delta U_{T_n}))^+)^2\Big]\nn\\
&=\E_t\Big[\sum_{n\in\mathbb{N}, \ t<T_n\leq T}((\delta Y_{j,T_n})^+)^2\Big],
\end{align}
where $U_t:=\int_0^t\int_{\cE}z N(\ds,\de)$, $\Delta U_{T_n}:=U_{T_n}-U_{T_n-}$, recalling that $\{T_n\}_{n\in\mathbb{N}}$ denotes of jump times of underlying Poisson point process.
Substituting \eqref{jumpterm} into \eqref{deltaY} yields,
\begin{align*}
((\delta Y_{i,t})^+)^2
&\leq c\mathbb{E}_t\int_t^T \sum_{i=1}^{\ell}((\delta Y_{i,s-})^+)^2\ds+\sum_{j\neq i}\E_t\Big[\sum_{n\in\mathbb{N}, \ t<T_n\leq T}((\delta Y_{j,T_n})^+)^2\Big].
\end{align*}
Since the jumps of $\delta Y$ are accountable, we can replace $\delta Y_{i,s-}$ by $\delta Y_{i,s}$ in the above integral to get
\begin{align}\label{maineq1}
((\delta Y_{i,t})^+)^2
&\leq c\mathbb{E}_t\int_t^T \sum_{i=1}^{\ell}((\delta Y_{i,s})^+)^2\ds+\sum_{j\neq i}\E_t\Big[\sum_{n\in\mathbb{N}, \ t<T_n\leq T}((\delta Y_{j,T_n})^+)^2\Big].
\end{align}

For any constant $h\in (0,T]$, set $$M(h):=\esssup\limits_{(t,i)\in[T-h,T]\times \cM}((\delta Y_{i,t})^+)^2,$$
which is finite since $\delta Y$ is bounded.
For any $t\in[T-h,T]$, we obtain from \eqref{maineq1} that
\begin{align*}
((\delta Y_{i,t})^+)^2
&\leq c\int_t^T \sum_{i=1}^{\ell}M(h)\ds+\sum_{j\neq i}\E_t\Big[\sum_{n\in\mathbb{N}, \ t<T_n\leq T}M(h)\Big]\\
&=c\ell M(h)(T-t)+M(h)\sum_{j\neq i}\E_t\int_t^T\int_{\cE}1N(\ds,\de)\\
&=c\ell M(h)(T-t)+M(h)(\ell-1)\nu(\cE)(T-t)\\
&\leq (c\ell+(\ell-1)\nu(\cE))M(h)h.
\end{align*}
Taking essential supreme over $(t,i)\in[T-h,T]\times \cM$ on both sides leads to
\begin{align}\label{Mineq}
M(h)&\leq (c\ell+(\ell-1)\nu(\cE))M(h)h.
\end{align}
Setting $h=\min\{1/(c\ell+(\ell-1)\nu(\cE)+1),T\}$ from now on. It then follows from above that $M(h)=0$, thus $\delta Y_{i,t}\leq0$ for all $t\in[T-h,T]$. Similarly, using $\delta Y_{i,T-h}\leq0$ and repeating the above argument on $[0\vee(T-2h),T-h]$, one can get $\delta Y_{i,t}\leq0$ for all $t\in[0\vee( T-2h),T-h]$.
Repeating this procedure, the desired comparison result follows.
\eof

\begin{remark}
If the inequalities in the conditions \ref{cond-boundary} and \ref{cond-size} are reversed, then so is the conclusion.
\end{remark}

\begin{remark}\label{remark:cond-gamma}
It is not hard to see that the condition \ref{cond-gamma} is equivalent to that
there exists a process $\gamma_{i}\in \linnu(0, T;\R)$ with $\gamma_{i}\geq-1$ such that
\begin{align*}
&\quad\;f_i(s,Y_{s-}, Z_{i,s}, \Phi_{1,s}, \cdots, \Phi_{i,s}, \cdots, \Phi_{\ell,s})\\
&\qquad\quad-f_i(s,Y_{s-}, Z_{i,s}, \Phi_{1,s}, \cdots, \overline\Phi_{i,s}, \cdots, \Phi_{\ell,s}) \\
&\leq \int_{\cE}\gamma_{i,s}(e)(\Phi_{i,s}(e)-\overline \Phi_{i,s}(e))\nu(\de).
\end{align*}
\end{remark}
Most of existing comparison theorems for BSDEJs require the condition $\gamma>-1$ or even stronger $\gamma>-1+\ep$ in order to utilize the Girsanov theorem; see, e.g., Barles, Buckdahn and Pardoux \cite{BBP} and Royer \cite{Royer}. Our requirement, namely $\gamma\geq -1$, is the same as Quenez and Sulem's \cite{QS}. But all these existing comparison theorems work for one-dimensional BSDEJs only.

\begin{remark} \label{remark:cond-growth}
The condition \ref{cond-growth} holds if,
for every $K>0$, there exists a constant $c>0$ (depending on $K$) such that
\begin{multline*}
f_{i}(s,y,z,\phi)-f_{i}(s, \overline y, \overline z, \overline \phi)\\
\qquad\leq c\Big(|y_{i}-\overline y_{i}| +\sum_{j\neq i} (y_{j}-\overline y_{j})^{+}+|z-\overline z|+\sum_{j\neq i}\int_{\cE}(y_{j}-\overline y_{j}+ \phi_{j}(e)-\overline \phi_{j}(e))^+\nu(\de)\Big)
\end{multline*}
holds for all $(y, z, \phi)$ and $(\overline y, \overline z, \overline \phi)\in \R^{\ell} \times \R^{n}\times\lnu(\R^{\ell})$ satisfying $ \phi_{i}\equiv\overline \phi_{i}$ and $|y|+|\overline y|\leq K$. Since $|y|+|\overline y|\leq K$, it is a locally Lipschitz condition w.r.t $y$.
The condition implies $f_{i}$ is increasing w.r.t $y_{j}$ and $\phi_{j}$ for every $j\neq i$. 
Also, the term $\sum_{j\neq i} (y_{j}-\overline y_{j})^{+}$ can be removed if there does exist jump, i.e. $\nu(\cE)>0$.
\end{remark}

\begin{remark} \label{remark:cond-growth-2}
We call a generator $f$ is Lipschitz in $(y,z,\phi)$ with Lipschitz constant $c$ if
\begin{align*}
|f(\omega,t,y,z,\phi)-f(\omega,t, \overline y, \overline z, \overline\phi)|\leq c(|y-\overline y|+|z-\overline z|+||\phi-\overline\phi||_{\nu})~\dpt
\end{align*}
holds for all $(y,z,\phi)$, $(\overline y,\overline z, \overline\phi)\in\R^\ell\times\R^{n\times \ell}\times L^{2,\nu}(\R^{\ell})$.
Then the condition \ref{cond-growth} holds if
\begin{enumerate}
\item $f_{i}(s,y,z,\phi)$ is Lipschitz in $(y,z,\phi)$;
\item $f_{i}(s,y,z,\phi)$ is increasing w.r.t $y_{j}$ for every $j\neq i$; and
\item there exists a constant $c>0$ such that
\begin{align*}
&f_{i}(s,Y_{s-},Z_{i,s},\Phi_{1,s},...,\Phi_{i-1,s},\overline\Phi_{i,s},\Phi_{i+1,s},...,\Phi_{\ell,s})\\
&\qquad-f_{i}(s,\overline Y_{1,s-},..., \overline Y_{i-1,s-}, Y_{i,s-},\overline Y_{i+1,s-},,...,\overline Y_{\ell,s-},Z_{i,s-},\overline \Phi_{s-})\\
&\leq c\sum_{j\neq i}\int_{\cE}(Y_{j,s-}+\Phi_{j,s}(e)-\overline Y_{j,s-}-\overline\Phi_{j,s}(e))^+\nu(\de).
\end{align*}
\end{enumerate}
\end{remark}

\begin{remark} In \eqref{deltaPhi} the condition $\gamma\in \linnu(0, T;\R^{\ell})$ can be replaced by
the following weaker one: there exist constants $0<h,\ep<1$ such that
$$\esssup_{t\in[0,T]}\mathbb{E}_t\int_t^{T\wedge (t+h)}\int_\cE |\gamma_s(e)|^{2} \nu(\de)\ds\leq 1-\ep.$$
This condition is satisfied, for instance, when $\int_\cE |\gamma_{\cdot}(e)|^{2} \nu(\de)\in L^{\infty}_{\mathbb{F}}(0, T;\R)$.
Indeed, the above condition implies, for $t\in[T-h,T]$,
\begin{align*}
\mathbb{E}_t\int_t^{T} ((\delta Y_{i,s-})^+)^2\int_\cE(\gamma_{i,s}(e)^{2}\vee 1)\nu(\de)\ds &\leq M(h)~\mathbb{E}_t\int_t^{T}\int_\cE(\gamma_{i,s}(e)^{2}+1)\nu(\de)\ds\\
&\leq M(h) (1-\ep+h\nu(\cE)) \leq (1-\ep/2)M(h),
\end{align*}
by choosing $h$ small enough.
This together with \eqref{deltaYsquare} and \eqref{deltaf} will lead to an estimate similar to \eqref{Mineq} in the above proof.
\end{remark}

\subsection{Comparison theorem for square integrable processes}
Theorem \ref{comparisonbound} requires the state processes to be bounded, which may be too restrictive for applications.
The following result relaxes this assumption to square integrable processes, but we have to in addition assume that both $f$ and $\overline f$ are globally Lipschitz.

\begin{theorem}
\label{comparison}
We shall use the same notations as in Theorem \ref{comparisonbound}.
Suppose, for all $ i\in \cM$,
\begin{align*}
(Y_i, Z_i, \Phi_i), (\overline Y_i, \overline Z_i, \overline\Phi_i)\in S^{2}_{\mathbb{F}}(0,T;\R)\times L^{2}_{\mathbb F}(0, T;\R^n)\times L^{2,\nu}_{\mathcal{P}}(0, T;\R),
\end{align*}
and they satisfy the BSDEJs \eqref{CY} and \eqref{CYoverline}.
Also suppose that, for every $i\in\cM$,
\begin{enumerate}
\item the conditions \ref{cond-boundary}, \ref{cond-gamma}, \ref{cond-growth} and \ref{cond-size} hold;
\item $f_{i}(\cdot, 0,0,0)$ and $\overline f_{i}(\cdot, 0,0,0)\in L^{2}_{\mathbb F}(0, T;\R)$;
\item both $f_{i}$ and $\overline f_i$ are Lipschitz in $(y,z,\phi)$.
\end{enumerate}
Then $Y_{i}\leq \overline Y_{i}$ for all $i\in\cM$.
\end{theorem}
\pf
For each $m\geq 1$ and $i\in\cM$, we denote
\begin{align*}
\xi^{m}_{i}=\xi_i\mathbf{1}_{|\xi|+|\overline \xi|\leq m}, ~~ f^{m}_{i}(t,y,z,\phi)=f_{i}(t,y,z,\phi)\mathbf{1}_{|f(t,0,0,0)|+|\overline f(t,0,0,0)|\leq m}, \\
\overline \xi^{m}_{i}=\overline \xi_i\mathbf{1}_{|\xi|+|\overline \xi|\leq m}, ~~ \overline f^{m}_{i}(t,y,z,\phi)=\overline f_{i}(t,y,z,\phi)\mathbf{1}_{|f(t,0,0,0)|+|\overline f(t,0,0,0)|\leq m}.
\end{align*}

Note that $\xi^{m}_{i}$, $\overline \xi^{m}_{i}$, $f^{m}_{i}(\cdot,0,0,0)$ and $\overline f^{m}_{i}(\cdot,0,0,0)$ are bounded by $m$ and the generators $f^m=(f^{m}_{1},...,f^{m}_{\ell})$ and $\overline f^m=(\overline f^{m}_{1},...,\overline f^{m}_{\ell})$ are both Lipschitz in $(y,z,\phi)$ with the same Lipschitz constant as $f$ and $\overline f$.
It then follows from \cite[Theorem 2.4]{TL} or \cite[Theorem 2.1, Proposition 2.2]{BBP} that the following BSDEJs:
\begin{align*}
Y_{i,t}^{m}&=\xi^{m}_{i}+\int_t^T f^{m}_{i}(s,Y_{s-}^{m}, Z_{i,s}^{m}, \Phi_{s}^{m})\ds\\
&\qquad\qquad\qquad-\int_t^T (Z_{i,s}^{m})^{\top}\dw_s-\int_t^T\int_{\cE}\Phi_{i,s}^{m}(e) \widetilde N(\ds,\de)~\dpt, \ i\in\cM,
\end{align*}
and
\begin{align*}
\overline Y_{i,t}^{m}&=\overline \xi^{m}_i+\int_t^T \overline f^{m}_i(s,\overline Y_{s-}^{m}, \overline Z_{i,s}^{m}, \overline \Phi_{s}^{m})\ds\\
&\qquad\qquad\qquad-\int_t^T (\overline Z_{i,s}^{m})^{\top}\dw_s-\int_t^T\int_{\cE}\overline \Phi_{i,s}^{m}(e) \widetilde N(\ds,\de)~\dpt, \ i\in\cM,
\end{align*}
admit unique solutions $(Y^{m}, Z^{m},\Phi^{m})$ and $(\overline Y^{m}, \overline Z^{m},\overline \Phi^{m})$ respectively, such that
\begin{align*}
(Y^{m}_i, Z^{m}_i,\Phi^{m}_i), ~(\overline Y^{m}_i, \overline Z^{m}_i,\overline \Phi^{m}_i)\in S^{2}_{\mathbb{F}}(0,T;\R)\times L^{2}_{\mathbb F}(0, T;\R^n)\times L_{\mathcal{P}}^{2,\nu}(0, T;\R)\ \mbox{ for all} \ i\in\cM.
\end{align*}

We temporally suppose that
\begin{align}\label{boundac}
Y^{m}_i,~~ \overline Y^{m}_i\in S^{\infty}_{\mathbb{F}}(0,T;\R)\ \mbox{ for all} \ i\in\cM.
\end{align}
Then applying Theorem \ref{comparisonbound} leads to
\begin{align}\label{Ybounded}
Y_{i}^{m}\leq\overline Y_{i}^{m} \ \mbox{for all} \ i\in\cM.
\end{align}
From \cite[Proposition 2.2]{BBP}, we know there is constant $c>0$ independent of $m$ such that
\begin{align*}
\E\Big[\sup_{0\leq t\leq T}|Y_t-Y_t^{m}|^2\Big]\leq c\E\Big[|\xi-\xi^m|^2+\int_0^T|f(t,Y_{t},Z_{t},\Phi_{t})-f^{m}(t,Y_{t},Z_{t},\Phi_{t})|^2\dt\Big],\\
\E\Big[\sup_{0\leq t\leq T}|\overline Y_t-\overline Y_t^{m}|^2\Big]\leq c\E\Big[|\overline \xi-\overline \xi^m|^2+\int_0^T|\overline f(t,\overline Y_{t},\overline Z_{t},\overline \Phi_{t})-\overline f^{m}(t,\overline Y_{t},\overline Z_{t},\overline \Phi_{t})|^2\dt\Big].
\end{align*}
These estimates together with the definitions of $\xi^m,\overline\xi^m,f^m,\overline f^m$ and the dominated convergence theorem lead to
\begin{align*}
\lim_{m\to \infty}\E\Big[\sup_{0\leq t\leq T}|Y_t-Y_t^{m}|^2
+ \sup_{0\leq t\leq T}|\overline Y_t-\overline Y_t^{m}|^2\Big]=0.
\end{align*}
Applying the elementary inequalities $(x^{+})^{2}\leq 2(y^{+})^{2}+2(x-y)^{2}$, $(x+y)^{2}\leq 2x^{2}+2y^{2}$ for $x$, $y\in\R$ and \eqref{Ybounded}, we have
\begin{align*}
&\quad\E\Big[\sup_{0\leq t\leq T}\sum_{i=1}^\ell[(Y_{i,t}-\overline Y_{i,t})^+]^2\Big]\\
&\leq \E\Big[2\sup_{0\leq t\leq T}\sum_{i=1}^\ell[(Y_{i,t}^m-\overline Y_{i,t}^{m})^+]^2 +2\sup_{0\leq t\leq T}\sum_{i=1}^\ell(Y_{i,t}-Y_{i,t}^{m}+\overline Y^m_{i,t}-\overline Y_{i,t})^2\Big]\\
&=\E\Big[2\sup_{0\leq t\leq T}\sum_{i=1}^\ell(Y_{i,t}-Y_{i,t}^{m}+\overline Y^m_{i,t}-\overline Y_{i,t})^2\Big]\\
&\leq \E\Big[4\sup_{0\leq t\leq T}\sum_{i=1}^\ell(Y_{i,t}-Y_{i,t}^{m})^2+4\sup_{0\leq t\leq T}\sum_{i=1}^\ell(\overline Y^m_{i,t}-\overline Y_{i,t})^2\Big].
\end{align*}
Sending $m\to \infty$ in the above, we get the desired result
$Y_{i}\leq \overline Y_{i}$ for all $i\in\cM$.

It remains to establish \eqref{boundac}. To this end,
let $\beta>0$ be a large constant to be chosen later. Applying It\^{o}'s formula to $e^{\beta t}(Y^m_{i,t})^2$, for each $i\in\cM$, yields
\begin{align*}
&\quad\; e^{\beta t}(Y^m_{i,t})^2+\E_t\int_t^T e^{\beta s}\Big(\beta (Y^m_{i,s})^2+|Z^m_{i,s}|^2+||\Phi^m_{i,s}||^2_{\nu}\Big)\ds\\
&=\E_t[e^{\beta T}(\xi^m_i)^2]+\E_t\int_t^T2e^{\beta s}Y^m_{i,s-} f^{m}_{i}(s,Y_{s-}^{m}, Z_{i,s}^{m}, \Phi_{s}^{m})\ds\\
&\leq m^{2}e^{\beta T}+\E_t\int_t^T 2e^{\beta s}|Y^m_{s-}| |f^m_{i}(s,Y_{s-}^{m}, Z_{i,s}^{m}, \Phi_{s}^{m})-f^m_i(s,0,0,0)| \ds\\
&\qquad\qquad\;+\E_t\int_t^T 2e^{\beta s}|Y^m_{s-}||f^m_i(s,0,0,0)|\ds\\
&\leq m^{2}e^{\beta T}+\E_t\int_t^T 2e^{\beta s}|Y^m_{s-}|c\Big(|Y_{s-}^{m}|+|Z_{i,s}^{m}|+||\Phi_{s}^{m}||_{\nu}\Big)\ds\\
&\qquad\qquad\;+\E_t\int_t^T e^{\beta s}|Y^m_{s-}|^{2}+\E_t\int_t^T e^{\beta s}|f^m_i(s,0,0,0)|^{2}\ds\\
&\leq m^2(1+ T)e^{\beta T}+\E_t\int_t^T e^{\beta s}\Big(c|Y_{s-}^{m}|^2+|Z_{i,s}^{m}|^{2}+\frac{1}{\ell}||\Phi_{s}^{m}||_{\nu}^{2}\Big)\ds,
\end{align*}
where the last constant $c$ does not depend on $t$, $\beta$ and $i$.
Canceling the common terms involving $|Z_{i,s}^{m}|^{2}$, we get
\begin{align*}
&\quad e^{\beta t}(Y^m_{i,t})^2+\E_t\int_t^T e^{\beta s}\Big(\beta(Y^m_{i,s})^2+||\Phi_{i,s}^{m}||_{\nu}^{2}\Big)\ds\\
&\leq m^2(1+ T)e^{\beta T}+\E_t\int_t^T e^{\beta s}\Big(c|Y_{s-}^{m}|^2+\frac{1}{\ell}||\Phi_{s}^{m}||_{\nu}^{2}\Big)\ds.
\end{align*}
Summing $i$ from $1$ to $\ell$ gives
\begin{align*}
&\quad e^{\beta t}|Y^m_{t}|^2+\E_t\int_t^T e^{\beta s}\Big(\beta |Y^m_{t}|^2+||\Phi^m_{t}||^2_{\nu}\Big)\ds\\
&\leq \ell m^2(1+ T)e^{\beta T}+\E_t\int_t^T e^{\beta s}\Big(c\ell|Y_{s-}^{m}|^2+||\Phi_{s}^{m}||_{\nu}\Big)\ds\\
&=\ell m^2(1+ T)e^{\beta T}+\E_t\int_t^T e^{\beta s}\Big(c\ell|Y_{s}^{m}|^2+||\Phi_{s}^{m}||_{\nu}\Big)\ds,
\end{align*}
where the last equation is due to the fact that the jumps of $Y$ are accountable.
By setting $\beta=c\ell$ and canceling the common integrals in the above estimate, we obtain $Y^{m}\in S^{\infty}_{\mathbb{F}}(0,T;\R^{\ell})$.
The assertion for $(\overline Y^{m}_i, \overline Z^{m}_i,\overline \Phi^{m}_i)$ in \eqref{boundac} can be similarly proved.
This completes the proof.
\eof

\section{A stochastic LQ control problem with jumps and the related two-dimensional BSDEJ}\label{section:LQ}
\subsection{Cone-constrained stochastic LQ control with jumps}
Consider the following $\R$-valued linear stochastic differential equation (SDE):
\begin{align}
\label{state}
\begin{cases}
\dd X_t=\left[A_tX_{t-}+B_t^{\top}u_t\right]\dt+\left[C_tX_{t-}+D_tu_t\right]^{\top}\dw_t\\
\qquad\qquad\qquad+\int_{\cE}\left[E_t(e)X_{t-}+F_t(e)^{\top}u_t\right] \widetilde N(\dt,\de), \ t\in[0,T], \\
X_0=x,
\end{cases}
\end{align}
where $A, \ B, \ C, \ D$ are all $\mathcal{P}$-measurable processes, and $E(\cdot), \ F(\cdot)$ are $ \mathcal{P}\otimes\mathcal{B}(\cE)$-measurable stochastic processes of suitable size, $x\in\R$ is known.

Let $\Pi$ be a given closed cone in $\R^m$, so if $u\in\Pi$, then $\lambda u\in\Pi$ for all $\lambda\geq 0$. It is used to represent the constraint set for controls.
The class of admissible controls is defined as the set
\begin{align*}
\mathcal{U}:= \Big\{u\in L^2_\mathbb{F}(0, T;\R^m)\;\Big|\; u_t \in\Pi,~\dpt\Big\}.
\end{align*}
If $u\in\mathcal{U}$, then \eqref{state} admits a unique solution $X$, and we call $(X, u)$ an admissible pair.

The cone-constrained stochastic LQ problem is stated as follows:
\begin{align}
\begin{cases}
\mathrm{Minimize} &\ J(x, u)\\
\mbox{subject to} &\ (X, u) \mbox{ is admissible for} \ \eqref{state},
\end{cases}
\label{LQ}%
\end{align}
where the cost functional is given as the following quadratic form
\begin{align}\label{costfunctional}
J(x, u):=&\mathbb{E}\Big[\int_0^T\Big(Q_tX_t^2+u_t^{\top}R_tu_t+2X_tS_t^{\top}u_t\Big)\dt+GX_T^2\Big].
\end{align}
The associated value function is defined as
\begin{align*}
V(x):=\inf_{u\in\mathcal{U}} J(x,u).
\end{align*}
Problem \eqref{LQ} is said to be solvable (at $x$), if there exists a control $u^*\in\mathcal{U}$ such that
\begin{align*}
-\infty<J(x, u^*)\leq J(x, u), \quad \forall\; u\in\mathcal{U},
\end{align*}
in which case, $u^*$ is called an optimal control for problem \eqref{LQ}, and the optimal value is
\begin{align*}
V(x)=J(x, u^*).
\end{align*}
Our aim is to solve problem \eqref{LQ}.

We put the following assumptions on the coefficients in this section.
\begin{assumption}[Bounded coefficients]
\label{assu1} It holds that
\begin{align*}
\begin{cases}
A\in L^{\infty}_{\mathbb{F}}(0, T;\R), \
B\in L_{\mathbb{F}}^\infty(0, T;\R^m), \
C \in L_{\mathbb{F}}^\infty(0, T;\R^n), \\
D\in L_{\mathbb{F}}^\infty(0, T;\R^{n\times m}), \
E\in\linnu(0, T;\R),\
F\in L_{\mathcal{P}}^{\infty,\nu}(0,T;\R^m),\\
Q\in L_{\mathbb{F}}^\infty(0, T;\R_+), \
R\in L_{\mathbb{F}}^\infty(0, T;\mathbb{S}^m),\
S\in L_{\mathbb{F}}^\infty(0, T;\R^m),\
G\in L_{\mathcal{F}_T}^\infty(\Omega;\R_+).
\end{cases}
\end{align*}
\end{assumption}
\begin{assumption} [Standard case]
\label{assu2}
It holds that
$\left(\begin{smallmatrix}
R & S\\
S^{\top} & Q
\end{smallmatrix}\right)\geq0$, and there exists a constant $\delta>0$ such that $R\geq \delta \mathbf{1}_{m}$, where $\mathbf{1}_m$ denotes the $m$-dimensional identity matrix.
\end{assumption}

\begin{assumption}[Singular case]
\label{assu3}
It holds that
$\left(\begin{smallmatrix}
R & S\\
S^{\top} & Q
\end{smallmatrix}\right)\geq0$ and there exists a constant $\delta>0$ such that $G\geq\delta$ and $D^{\top}D+\int_{\cE}F(e)F(e)^{\top}\nu(\de)\geq\delta \mathbf{1}_m$.
\end{assumption}

\subsection{Coupled SRE with jumps}
Nowadays, it is well known that solutions to stochastic LQ problems depend heavily on the solvability of the related SREs. We now introduce the associated SRE for our problem \eqref{LQ}.\footnote{We will give a heuristic derivation in Appendix \ref{appnA} for the readers' convenience. See also Dong \cite{Dong} for a special SRE with single jump stems from the theory of filtration enlargement.}

For $(\omega,t,v,P_i,\Lambda,\Gamma_i)\in \Omega\times[0,T]\times\Pi\times \R\times\R^m\times L^{\infty,\nu}(\R)$, $i=1,2$ define the following mappings:
\begin{align*}
H_1(\omega,t,v,P_1,P_2,\Lambda,\Gamma_1, \Gamma_2)&:=v^{\top}(R+P_1D^{\top}D)v+2(P_1(B+D^{\top}C)+D^{\top}\Lambda+S)^{\top}v\\
&\qquad+\int_{\cE}\Big[(P_1+\Gamma_1)\Big(((1+E+F^{\top}v)^+)^2-1\Big)-2P_1(E+F^{\top}v)\\
&\qquad\qquad\quad+(P_2+\Gamma_2) ((1+E+F^{\top}v)^-)^2\Big]\nu(\de), \\
H_2(\omega,t,v,P_1,P_2,\Lambda,\Gamma_1, \Gamma_2)&:=v^{\top}(R+P_2D^{\top}D)v-2(P_2(B+D^{\top}C)+D^{\top}\Lambda+S)^{\top}v\\
&\qquad+\int_{\cE}\Big[(P_2+\Gamma_2)\Big(((-1-E+F^{\top}v)^-)^2-1\Big)+2P_2(-E+F^{\top}v)\\
&\qquad\qquad\quad+(P_1+\Gamma_1)\Gamma((-1-E+F^{\top}v)^+)^2\Big]\nu(\de),
\end{align*}
and set
\begin{align}\label{Hstar1}
H_1^{*}(\omega,t,P_1,P_2,\Lambda,\Gamma_1, \Gamma_2):=\inf_{v\in\Pi}
H_1(\omega,t,v,P_1,P_2,\Lambda,\Gamma_1, \Gamma_2),\\
H_2^{*}(\omega,t,P_1,P_2,\Lambda,\Gamma_1, \Gamma_2):=\inf_{v\in\Pi}
H_2(\omega,t,v,P_1,P_2,\Lambda,\Gamma_1, \Gamma_2).\label{Hstar2}
\end{align}
The associated SRE for our problem \eqref{LQ} is given as follows:
\begin{align}\label{P1}
\begin{cases}
\dd P_{1,t}=-\Big[(2A+C^{\top}C)P_{1,t-}+2C^{\top}\Lambda_{1,t}+Q+
H_1^*(t,P_{1,t-},P_{2,t-},\Lambda_{1,t},\Gamma_{1,t},\Gamma_{2,t})\Big]\dt\\
\hfill+\Lambda_{1,t}^{\top}\dw+\int_{\cE}\Gamma_{1,t}(e) \widetilde N(\dt,\de),\qquad\quad\\
\dd P_{2,t}=-\Big[(2A+C^{\top}C)P_{2,t-}+2C^{\top}\Lambda_{2,t}+Q+
H_2^*(t,P_{1,t-},P_{2,t-},\Lambda_{2,t},\Gamma_{1,t},\Gamma_{2,t})\Big]\dt\\
\hfill+\Lambda_{2,t}^{\top}\dw+\int_{\cE}\Gamma_{2,t}(e) \widetilde N(\dt,\de),\qquad\quad\\
P_{1,T}=G, \ P_{2,T}=G,\\
P_{1,t}\geq 0, \ P_{1,t-}+\Gamma_{1,t}\geq 0, \ P_{2,t}\geq 0, \ P_{2,t-}+\Gamma_{2,t}\geq 0.
\end{cases}
\end{align}
This is a new two-dimensional coupled nonlinear BSDEJ.

\begin{remark}
Hu and Zhou \cite{HZ} studied a cone-constrained LQ problem without jumps; the associated SREs \cite[Eq. (3.5) and (3.6)]{HZ} are decoupled, so that one can solve $P_1$ and $P_2$ separately. As is well-known $P_1$ and $P_2$ correspond to the optimal value with positive and negative initial state. When there is no jump in the model, the optimal state process does not change sign, so that only one of $P_1$ and $P_2$ is involved.
Therefore, they are decoupled.

Things become notably different in models with jumps. Because of jumps, the sign of the optimal state process may switch between positive and negative values, so $P_1$ and $P_2$ are coupled together and one cannot treat them separately.
So our SRE \eqref{P1} is actually a system of coupled BSDEJs whose solvability is far from trivial compared to the decoupled BSDEJs in \cite[Eq. (3.5) and (3.6)]{HZ}.

If all the coefficients in Assumption \ref{assu1} are predictable with respect to the Brownian filtration, then $\Gamma_1=\Gamma_2=0$ and the SRE becomes a two-dimensional coupled BSDE without jumps. Even without jumps, the BSDE is still new and cannot be covered by existing results on multi-dimensional BSDEs; see, e.g., Fan, Hu and Tang \cite{FHT}, Hu and Tang \cite{HT}.
\end{remark}

\begin{remark}
If $\Pi$ is symmetric, namely, $-v\in\Pi$ whenever $v\in\Pi$, then $H_1^*=H_2^*$ and \eqref{P1} will degenerate to one equation since $(P_1,\Lambda_1,\Gamma_1)=(P_2,\Lambda_2,\Gamma_2)$. In particular, if there is no control constraint, that is, $\Pi=\R^m$, then both $H_1^*$ and $H_2^*$ are equal to
\begin{align*}
(P+\Gamma)E^2+2\Gamma E&+\Big(P(B+D^{\top}C)+D^{\top}\Lambda+S+\int_{\cE}((P+\Gamma)E+\Gamma)F\nu(\de)\Big)^{\top}\\
&\qquad\times \Big(R+PD^{\top}D+\int_{\cE}(P+\Gamma)FF^{\top}\nu(\de)\Big)^{-1}\\
&\qquad\times\Big(P(B+D^{\top}C)+D^{\top}\Lambda+S+\int_{\cE}((P+\Gamma)E+\Gamma)F\nu(\de)\Big).
\end{align*}
Under $\Pi=\R^m$, Zhang, Dong and Meng \cite{ZDM} addressed the solvability of the matrix-valued SREJ under the assumption $R\geq \delta\mathbf{1}_m$ and $S\equiv 0$. By contrast, we will solve the BSDEJ \eqref{P1} in both standard and singular cases for general cone $\Pi$.
\end{remark}

\begin{definition}\label{def}
A stochastic process $(P_1,\Lambda_{1},\Gamma_1,P_2,\Lambda_{2},\Gamma_2)$ is called a solution to the BSDEJ \eqref{P1} if it satisfies \eqref{P1}, and $(P_i,\Lambda_{i},\Gamma_i)\in S^{\infty}_{\mathbb{F}}(0,T;\R)\times L^{2}_{\mathbb{F}}(0,T;\R^m)\times \linnu(0, T;\R)$, $i=1,2$. The solution is called nonnegative if $P_{i}\geq0$, and called uniformly positive if $P_{i}\geq c$ for some deterministic constant $c>0$, $i=1,2$.
\end{definition}

\subsection{Existence of solution to the BSDEJ \eqref{P1}}
Dong \cite{Dong} constructed a solution to a SRE with single jump using two recursive systems of BSDEs driven only by Brownian motions. His decomposition approach is tailor made in the filtration enlargement framework, hence fails in the Poisson random measure model which accommodates probably accountable jumps.

Czichowsky and Schweizer \cite{CS} characterized the optimal value process of a cone-constrained mean-variance problem in terms of a coupled system of BSDEs \cite[Eq.(4.18)]{CS} in a semimartingale model. They claimed in \cite[Remark 4.8]{CS} that \emph{``Due to the coupling term coming from $\mathfrak{h}$, the BSDE system
(4.18) is very complicated. It has a nonlinear non-Lipschitz generator plus a
generator with jumps, so that finding a solution by general BSDE techniques seems a
formidable challenge''}. We now respond to this formidable challenge in the Wiener-Poisson world by providing a proof of the existence of solution to \eqref{P1} by pure BSDE techniques.
\begin{theorem}[Existence in Standard case]
\label{Th:standard}
Suppose Assumptions \ref{assu1}, \ref{assu2} hold, then the BSDEJ \eqref{P1} admits a nonnegative solution $(P_1,\Lambda_1,\Gamma_1,P_2,\Lambda_2,\Gamma_2)$.
\end{theorem}
\pf
For $k=1,2,...$, define maps
\begin{align}\label{defHk}
H_1^{k}(\omega,t,P_{1},P_{2},\Lambda_{1},\Gamma_{1},\Gamma_2):=\inf_{v\in\Pi,|v|\leq k}
H_1(\omega,t,v,P_{1},P_{2},\Lambda_{1},\Gamma_{1},\Gamma_2),\\
H_2^{k}(\omega,t,P_{1},P_{2},\Lambda_{2},\Gamma_{1},\Gamma_2):=\inf_{v\in\Pi,|v|\leq k}
H_2(\omega,t,v,P_{1},P_{2},\Lambda_{2},\Gamma_{1},\Gamma_2).
\end{align}
Then they are uniformly Lipschitz in $(P_{1},\Lambda_{1},\Gamma_1,P_{2}, \Lambda_2,\Gamma_2)$ and decreasingly approach to\\
$H_1^{*}(\omega,t,P_{1},\Lambda_{1},\Gamma_1,P_{2},\Gamma_2)$ and $H_2^{*}(\omega,t,P_{2},\Lambda_{2},\Gamma_{2},P_{1},\Gamma_1)$ respectively as $k$ goes to infinity.

For each $k$,
the following BSDE
\begin{align}\label{Ptrun}
\begin{cases}
\dd P_{1,t}^k=-\Big[(2A+C^{\top}C)P_{1,t-}^k+2C^{\top}\Lambda_{1,t}^k+Q
+H_1^{k}(t,P_{1,t-}^k,P_{2,t-}^k,\Lambda_{1,t}^k,\Gamma_{1,t}^k,\Gamma_{2,t}^k)\Big]\dt\\
\hfill+(\Lambda_{1,t}^k)^\top \dw+\int_{\cE}\Gamma_{1,t}^k(e) \widetilde N(\dt,\de),\qquad\quad\\
\dd P_{2,t}^k=-\Big[(2A+C^{\top}C)P_{2,t-}^k+2C^{\top}\Lambda_{2,t}^k+Q
+H_2^{k}(t,P_{1,t-}^k,P_{2,t-}^k,\Lambda_{2,t}^k,\Gamma_{1,t}^k,\Gamma_{2,t}^k)\Big]\dt\\
\hfill+(\Lambda_{2,t}^k)^\top \dw+\int_{\cE}\Gamma_{2,t}^k(e) \widetilde N(\dt,\de),\qquad\quad\\
P_{1,T}^k=G, \ P_{2,T}^k=G,
\end{cases}
\end{align}
is a two-dimensional BSDEJ with a Lipschitz generator, so by \cite[Lemma 2.4]{TL}, it admits a unique solution $(P_1^k,\Lambda_1^k,\Gamma_1^k,P_2^k,\Lambda_2^k,\Gamma_2^k)$ such that $$(P_i^k,\Lambda_i^k,\Gamma_i^k)\in S^{2}_{\mathbb{F}}(0, T;\R)\times L^{2}_{\mathbb{F}}(0, T;\R^n)\times \ltwonu(0, T;\R), \ i=1,2.$$
From the definition of $H_1^{k}$, we have
\begin{align*}
&\qquad H_1^{k}(t,P_{1},P_{2},\Lambda_{1},\Gamma_{1},\Gamma_{2})
-H_1^{k}(t,P_{1},P'_{2},\Lambda_{1},\Gamma_{1},\Gamma'_{2})\\
&\leq\sup_{v\in\Pi,|v|\leq k}\int_{\cE} ((1+E+F^{\top}v)^-)^2(P_2+\Gamma_2(e)-P'_{2}-\Gamma'_{2}(e))\nu(\de)\\
&\leq c_{k}\int_{\cE} (P_2+\Gamma_2(e)-P'_{2}-\Gamma'_{2}(e))^+\nu(\de),
\end{align*}
and
\begin{align*}
&\qquad H_1^{k}(t,P_{1},P_{2},\Lambda_{1},\Gamma_{1},\Gamma_2)
-H_1^{k}(t,P_{1},P_{2},\Lambda_{1},\Gamma'_{1},\Gamma_2)\\
&\leq \sup_{v\in\Pi,|v|\leq k}\int_{\cE}(\Gamma_{1}(e)-\Gamma'_{1}(e))\Big(((1+E+F^{\top}v)^+)^2-1\Big)\nu(\de)\\
&\leq \sup_{v\in\Pi,|v|\leq k}\int_{\cE}(\Gamma_{1}(e)-\Gamma'_{1}(e))((1+E+F^{\top}v)^+)^2\nu(\de)+ \int_{\cE}|\Gamma_{1}(e)-\Gamma'_{1}(e)|\nu(\de)\\
&\leq c_{k} \int_{\cE}(\Gamma_{1}(e)-\Gamma'_{1}(e))^{+}\nu(\de)+ \int_{\cE}|\Gamma_{1}(e)-\Gamma'_{1}(e)|\nu(\de),
\end{align*}
where $c_{k}<\infty$ is defined as
\begin{align*}
c_k&=\esssup_{v\in\Pi,|v|\leq k}|1+E+F^{\top}v|^2,~\nu(\de)\textrm{-a.e..}
\end{align*}
Similar estimates for $H_2^{k}$ can be established.
Hence according to
Theorem \ref{comparison}, $P^k_i$ is decreasing in $k$, for $i=1,2$.

Next, we show that the sequence $\{P_i^k\}_{k=1,2,...}$ is nonnegative and uniformly bounded from above, for $i=1,2$.

From Assumption \ref{assu1}, there exists a constant $c>0$ such that
\begin{align*}
2A+C^{\top}C+\int_{\cE}E(e)^2\nu(\de)\leq c, \ Q\leq c, \ G\leq c.
\end{align*}
It is easy to check that $(\overline P_{1,t},\overline\Lambda_{1,t},\overline\Gamma_{1,t})=(\overline P_{2,t},\overline\Lambda_{2,t},\overline\Gamma_{2,t})=((c+1)e^{c(T-t)}-1,0,0)$ satisfies
the following two-dimensional BSDEJ
\begin{align}
\label{Povli}
\begin{cases}
\dd\overline P_1=-\Big[c\overline P_1+C^{\top}\Lambda_1+c+\int_{\cE}\Big(\overline \Gamma_1\big(((1+E)^+)^2)-1\big)+\overline \Gamma_2((1+E)^-)^2\Big)\nu(\de)\Big]\dt\\
\hfill+\overline \Lambda_1^\top \dw+\int_{\cE}\overline \Gamma_1(e) \widetilde N(\dt,\de),\qquad\quad\\
\dd\overline P_2=-\Big[c\overline P_2+C^{\top}\Lambda_2+c+\int_{\cE}\Big(\overline \Gamma_2\big(((1+E)^+)^2)-1\big)+\overline \Gamma_2((1+E)^-)^2\Big)\nu(\de)\Big]\dt\\
\hfill+\overline \Lambda_2^\top \dw+\int_{\cE}\overline \Gamma_2(e) \widetilde N(\dt,\de),\qquad\quad\\
\overline P_{1,T}=c,~\overline P_{2,T}=c.
\end{cases}
\end{align}
By the definition of $H_1^{k}$, we have
\begin{align*}
H_1^{k}(t,\overline P_{1},\overline P_{2},\overline\Lambda_{1},\overline\Gamma_{1},\overline\Gamma_{2}) &\leq H_1(t,0,\overline P_{1},\overline P_{2},\overline\Lambda_{1},\overline\Gamma_{1},\overline\Gamma_{2})\\
&=\int_{\cE}\Big(\overline P_1E^2+\overline \Gamma_1\big(((1+E)^+)^2)-1\big)+\overline \Gamma_1((1+E)^-)^2\Big)\nu(\de),
\end{align*}
so
\begin{align*}
&\quad (2A+C^{\top}C)\overline P_{1}+2C^{\top}\overline\Lambda_{1}+Q
+H_1^{k}(t,\overline P_{1},\overline P_{2},\overline\Lambda_{1},\overline\Gamma_{1},\overline\Gamma_{2})\\
&\leq c\overline P_1+C^{\top}\overline\Lambda_1+c+\int_{\cE}\Big(\overline \Gamma_1\big(((1+E)^+)^2)-1\big)+\overline \Gamma_1((1+E)^-)^2\Big)\nu(\de).
\end{align*}
Similarly, we have
\begin{align*}
&\quad (2A+C^{\top}C)\overline P_{2}+2C^{\top}\overline\Lambda_{2}+Q
+H_2^{k}(t,\overline P_{1},\overline P_{2},\overline\Lambda_{2},\overline\Gamma_{1},\overline\Gamma_{2})\\
&\leq c\overline P_2+C^{\top}\Lambda_2+c+\int_{\cE}\Big(\overline \Gamma_2\big(((1+E)^+)^2)-1\big)+\overline \Gamma_2((1+E)^-)^2\Big)\nu(\de).
\end{align*}
Keeping the above two inequalities in mind, applying Theorem \ref{comparison} to BSDEJs \eqref{Ptrun} and \eqref{Povli}, we have for $i=1,2$, $k=1,2,...$
\begin{align}
\label{upperbound}
P_{i,t}^k\leq \overline P_{i,t}\leq M,
\end{align}
where $M:=(c+1)e^{cT}-1$.

On the other hand, notice that $(\underline P_{1,t},\underline\Lambda_{1,t},\underline\Gamma_{1,t})=(\underline P_{2,t},\underline\Lambda_{2,t},\underline\Gamma_{2,t}):=(0,0,0)$ satisfies
\begin{align*}
\begin{cases}
\dd\underline P=\underline \Lambda^{\top}\dw+\int_{\cE}\underline\Gamma(e) \widetilde N(\dt,\de),\\
\underline P_T=0,
\end{cases}
\end{align*}
and
\begin{align*}
&\quad (2A+C^{\top}C)\underline P_{1}+2C^{\top}\underline\Lambda_{1}+Q
+H_1^{k}(t,\underline P_{1},\underline P_{2},\underline\Lambda_{1},\underline\Gamma_{1},\underline\Gamma_{2})\\
&\geq Q+\inf_{v\in\R^m}(v^{\top}Rv+2S^{\top}v)=Q-S^{\top}R^{-1}S\geq0,
\end{align*}
thanks to Assumption \ref{assu2}. Hence, by Theorem \ref{comparison} again,
\begin{align}\label{lowerbound}
P^{k}_{i,t}\geq \underline P_t=0, \ i=1,2,~k=1,2,...
\end{align}

Notice, for $i=1,2$,
\begin{align*}
\E\int_0^T\int_{\cE}\mathbf{1}_{\{P_{i,t-}^k+\Gamma_{i,t}^k(e)<0\}} \nu(\de)\dt&
=\E\int_0^T\int_{\cE}\mathbf{1}_{\{P_{i,t-}^k+\Gamma_{i,t}^k(e)<0\}}N(\dt,\de)\\
&=\E\Big[\sum_{n\in\mathbb{N},T_n\leq T}\mathbf{1}_{\{P_{i,T_n-}^k+\Gamma_{i,T_n}^k(-\Delta U_{T_n})<0\}}\Big]\\
&=\E\Big[\sum_{n\in\mathbb{N},T_n\leq T}\mathbf{1}_{\{P_{i,T_n}^k<0\}}\Big]=0,
\end{align*}
where $U_t=\int_0^t\int_{\cE}eN(\dt,\de)$ and $\Delta U_{T_n}= U_{T_n}- U_{T_n-}$,
hence,
\begin{align}\label{PGammageq}
P_{i,t-}^k+\Gamma_{i,t}^k\geq0.
\end{align}
Similarly, we can establish
\begin{align}\label{PGammaleq}
P_{i,t-}^k+\Gamma_{i,t}^k\leq M.
\end{align}
Now we obtain
$$-M\leq - P_{i,t-}^k\leq \Gamma_{i,t}^k\leq M- P_{i,t-}^k\leq M.
$$
Hence, $\Gamma_i^k$, $k=1,2,\cdots,$ are uniformly bounded by $M$, and thus belong to $\linnu(0, T;\R) $.

Since $P_{i}^k$ is decreasing w.r.t. $k$, we can define $P_{i,t}:=\lim_{k\to \infty}P_{i,t}^k, \ i=1,2$. Combining \eqref{upperbound} and \eqref{lowerbound}, it follows
\[
0\leq P_{i,t}\leq M, \ i=1,2,~ t\in[0,T].
\]
Applying It\^{o}'s formula to $(P_{1,t}^k)^2$, we deduce that
\begin{align*}
\begin{cases}
\dd\;(P_{1,t}^k)^2=\Big\{-2P_{1}^{k}\Big[(2A+C^{\top}C)P_{1,t-}^k+2C^{\top}\Lambda_{1,t}^k+Q
+H_1^{k}(t,P_{1,t-}^k,P_{2,t-}^k,\Lambda_{1,t}^k,\Gamma_{1,t}^k,\Gamma_{2,t}^k)\Big]\\
\qquad\qquad\qquad
+|\Lambda_{1}^k|^2+\int_{\cE}\Gamma_{1}^k(e)^2 \nu(\de)\Big\}\dt\\
\hfill+2P_{1}^k(\Lambda_{1}^k)^{\top}\dw+\int_{\cE}[(P_{1,t-}^k+\Gamma_{1,t}^k(e))^2-(P_{1,t-}^k)^2] \widetilde N(\dt,\de),\qquad\quad\\
(P_{1,T}^k)^2=G^2.
\end{cases}
\end{align*}
Since $0\leq P_{i}^k, P_{i}^k+\Gamma_i^k\leq M, i=1,2$, and $$H_1^{k}\leq \int_{\cE}\Big[(P_{1}+\Gamma_{1})\Big(((1+E)^+)^2-1\Big)-2P_{1}E+(P_{2}+\Gamma_{2})((1+E)^-)^2\Big]\nu(\de)\leq c,$$
by taking expectation on both sides in above and integrating over $[0,T]$, we have
\begin{align}\label{L2bound}
&\quad (P_{1,0}^k)^2+\frac{1}{2}\E\int_0^T|\Lambda_{1}^k|^2\ds+\E\int_0^T\int_{\cE}\Gamma_{1}^k(e)^2\nu(\de)\ds \leq c,
\end{align}
where $c>0$ is a constant independent of $k$.
Therefore, the sequence $(\Lambda_{1}^k, \Gamma_{1}^k)$, $k=1,2,\cdots,$ is bounded in $L^2_{\mathbb{F}}(0,T;\R^n)\times \ltwonu(0, T;\R)$, thus we can extract a subsequence (which is still denoted by $(\Lambda_{1}^k, \Gamma_{1}^k)$) converging in the weak sense to some $(\Lambda_{1}, \Gamma_{1})\in L^2_{\mathbb{F}}(0,T;\R^n)\times \ltwonu(0, T;\R)$.
Similar considerations applying to $(P_{2,t}^k)^2$ yield some $(\Lambda_{2}, \Gamma_{2})\in L^2_{\mathbb{F}}(0,T;\R^n)\times \ltwonu(0, T;\R)$ which is the weak limit of $(\Lambda_{2}^k, \Gamma_{2}^k)$.

Following Kobylanski's argument \cite[Proposition 2.4]{Ko}
(see also Antonelli and Mancini \cite[Theorem 1]{AM}, Kohlmann and Tang \cite[Theorem 2.1]{KT}), we establish in Appendix \ref{appnB} the following strong convergence:
\begin{lemma}\label{thlemma1}
It holds that
\begin{align}\label{strongcon}
\lim_{k\to \infty}\E\int_0^T|\Lambda_i^k-\Lambda_i|^2\dt=0, \ \lim_{k\to \infty}\E\int_0^T\int_{\cE}|\Gamma_i^k-\Gamma_i|^2\nu(\de)\dt=0, \ i=1,2.
\end{align}
Furthermore, $(P_1,\Lambda_1,\Gamma_1,P_2,\Lambda_2,\Gamma_2)$ is a nonnegative solution to the BSDEJ \eqref{P1}.
\end{lemma}
This completes the proof.
\eof

\begin{theorem}[Existence in Singular case]\label{Th:singular}
Suppose Assumptions \ref{assu1}, \ref{assu3} hold, then the BSDEJ \eqref{P1} admits a uniformly positive solution $(P_1,\Lambda_1,\Gamma_1,P_2,\Lambda_2,\Gamma_2)$.
\end{theorem}
\pf
Similar to the proof of Theorem \ref{Th:standard}, one can show the existence of a nonnegative solution $(P_1,\Lambda_1,\Gamma_1,P_2, \Lambda_2, \Gamma_2)$ to the BSDEJ \eqref{P1}, so we omit the details.
We only give a sketch on how to find a uniformly positive lower bound for such a solution.

Under Assumptions \ref{assu1}, \ref{assu3}, there exists constant $c_2>0$, such that
\[
2A+C^{\top}C+\int_{\cE}E^2\nu(\de)-\delta^{-1}\Big|B+D^{\top}C\pm\int_{\cE}EF\nu(\de)\Big|^2\geq -c_2,
\]
where $\delta$ is the constant in Assumption \ref{assu3}.
Notice $(\underline P_{1,t},\underline\Lambda_{1,t},\underline\Gamma_{1,t})=(\underline P_{2,t},\underline\Lambda_{2,t},\underline\Gamma_{2,t}):=(\delta e^{-c_2(T-t)},0,0)$ solves the following BSDEJ
\begin{align}\label{Plowersingu}
\begin{cases}
\dd\underline P=-(-c_2\underline P+C^{\top}\underline\Lambda)\dt+\underline\Lambda^{\top}\dw+\int_{\cE}\underline\Gamma(e) \widetilde N(\dt,\de),\\
\underline P_{T}=\delta.
\end{cases}
\end{align}
We have the following estimates,
\begin{align*}
&\qquad H_1^{k}(t,\underline P_{1},\underline P_{2},\underline \Lambda_{1},\underline \Gamma_1,\underline \Gamma_2)\\
&\geq\inf_{v\in\R^m}
H_1(t,v,\underline P_{1},\underline P_{2},\underline \Lambda_{1},\underline \Gamma_1,\underline \Gamma_2)\\
&=\inf_{v\in\R^m}
H_1(t,v,\underline P_{1},\underline P_{2},0,0,0)\\
&\geq \inf_{v\in\R^m}\Big[v^{\top}Rv+2S^{\top}v\Big]+\underline P_1\int_{\cE}E^2\nu(\de)\\
&\qquad+
\underline P_{1}\inf_{v\in\R^m}\Big[v^{\top}(D^{\top}D+\int_{\cE}F(e)F(e)^{\top}\nu(\de))v+2\Big(B+D^{\top}C+\int_{\cE}EF\nu(\de)\Big)^{\top}v\Big]\\
&\geq -Q+\underline P_1\Big[\int_{\cE}E^2\nu(\de)-\delta^{-1}\Big|B+D^{\top}C+\int_{\cE}EF\nu(\de)\Big|^2\Big],
\end{align*}
where we used $\underline\Lambda_1=0$, $\underline\Gamma_1=\underline\Gamma_2=0$ in the equality, $\underline P_1=\underline P_2>0$ in the second inequality, and
$$\left(\begin{smallmatrix}
R & S\\
S^{\top} & Q
\end{smallmatrix}\right)\geq0,\quad D^{\top}D+\int_{\cE}F(e)F(e)^{\top}\nu(\de)\geq\delta \mathbf{1}_m$$ in the last inequality. Similar result also holds for $H_2^{k}(t,\underline P_{1},\underline P_{2},\underline \Lambda_2,\underline \Gamma_1,\underline \Gamma_{2})$.

Applying Theorem \ref{comparison} to the BSDEJs \eqref{Ptrun} and \eqref{Plowersingu}, we get, for $i=1,2$,
\begin{align}
P_{i,t}^k\geq \underline P_{i,t}=\delta e^{-c_2(T-t)}\geq \delta e^{-c_2T}, \ t\in[0,T],
\end{align}
which leads to the desired lower bound.
\eof

\subsection{Solution to the LQ problem \eqref{LQ}}
In this subsection we will present an explicit solution to the LQ problem \eqref{LQ} in terms of solutions to the BSDEJ \eqref{P1}.

For $P_i>0$, $\Lambda_i\in\R^n$, $\Gamma_i\in L^{2,\nu}, \ i=1,2$, define
\begin{align}\label{hatv}
&\hat v_1(\omega,t,P_{1},P_{2},\Lambda_{1},\Gamma_{1},\Gamma_2)
=\argmin_{v\in\Pi}H_1(\omega,t,v,P_{1},P_{2},\Lambda_{1},\Gamma_{1},\Gamma_2),\nn\\
&\hat v_2(\omega,t,P_{1},P_{2},\Lambda_{2},\Gamma_{1},\Gamma_2)
=\argmin_{v\in\Pi}H_1(\omega,t,v,P_{1},P_{2},\Lambda_{2},\Gamma_{1},\Gamma_2).
\end{align}

\begin{theorem}\label{Th:verif}
Let $(P_i,\Lambda_i,\Gamma_i)\in S^{\infty}_{\mathbb{F}}(0,T;\R)\times L^{2}_{\mathbb{F}}(0,T;\R^m) \times \linnu(0, T;\R), \ i=1,2$, be a nonnegative (in Standard case), or uniformly positive (in Singular case) solution to the BSDEJ \eqref{P1}. Then the state feedback control given by
\begin{align*}
u^*(t,X)=\hat v_1(\omega,t,P_{1,t-},P_{2,t-},\Lambda_{1,t},\Gamma_{1,t},\Gamma_{2,t})X_{t-}^++\hat v_2(\omega,t,P_{1,t-},P_{2,t-},\Lambda_{1,t},\Gamma_{1,t},\Gamma_{2,t})X_{t-}^-,
\end{align*}
is optimal for the LQ problem \eqref{LQ}. Moreover, the optimal value is
\begin{align*}
V(x)=P_{1,0}(x^+)^2+P_{2,0}(x^-)^2.
\end{align*}
\end{theorem}
The proof of Theorem \ref{Th:verif} is standard, and thus omitted here; please see \cite[Theorem 5.2]{HZ} or \cite[Theorem 5.2]{ZDM} for the standard verification argument.

As a byproduct of Theorem \ref{Th:verif}, we have the following uniqueness result.
\begin{theorem}\label{uniqueness}
Suppose Assumptions \ref{assu1} and \ref{assu2} (resp. Assumptions \ref{assu1} and \ref{assu3}) hold, then the BSDEJ \eqref{P1} admits at most one nonnegative (resp. uniformly positive) solution.
\end{theorem}
It seems a challenging task to establish this result by pure BSDE techniques.

\section{Application to mean-variance portfolio selection problem} \label{section:MV}
Consider a financial market consisting of a risk-free asset (the money market
instrument or bond) whose price is $S_{0}$ and $m$ risky securities (the
stocks) whose prices are $S_{1}, \ldots, S_{m}$. And assume $m\leq n$, i.e. the number of risky securities is no more than the dimension of the Brownian motion.
The asset prices $S_k$, $k=0, 1, \ldots, m, $ are driven by stochastic differential equations (SDEs):
\begin{align*}
\begin{cases}
\dd S_{0,t}=r_tS_{0,t}\dt, \\
S_{0,0}=s_0,
\end{cases}
\end{align*}
and
\begin{align*}
\begin{cases}
\dd S_{k,t}=S_{k,t}\Big((\mu_{k,t}+r_t)\dt+\sum\limits_{j=1}^n\sigma_{kj,t}\dw_{j,t}+\int_{\cE}F_{k,t}(e) \widetilde N(\dt,\de)\Big), \\
S_{k,0}=s_k,
\end{cases}
\end{align*}
where, for every $k=1, \ldots, m$, $r$ is the interest rate process, $\mu_k$, $\sigma_k:=(\sigma_{k1}, \ldots, \sigma_{kn})$ and $F_k$ are the mean excess return rate process and volatility rate process of the $k$-th risky security.

Define the vectors
$\mu=(\mu_1, \ldots, \mu_m)^{\top}$, $F=(F_1, \ldots, F_m)^{\top}$
and matrix
\begin{align*}
\sigma=
\left(
\begin{array}{c}
\sigma_1\\
\vdots\\
\sigma_m\\
\end{array}
\right)
\equiv (\sigma_{kj})_{m\times n}, \ \text{for}\ \text{each} \ i\in\cM.
\end{align*}
We shall assume, in this section,
\begin{assumption} \label{assu4}
The interest rate $r$ is a bounded deterministic measurable function of $t$,
$$\mu\in L_{\mathbb{F}}^\infty(0, T;\R^m), \
\sigma\in L_{\mathbb{F}}^\infty(0, T;\R^{m\times n}), \
F\in L_{\mathcal{P}}^{\nu,\infty}(0,T;\R^m),$$
and there exists a constant $\delta>0$ such that
$\sigma\sigma^{\top}+\int_{\cE}F(e)F(e)^{\top}\nu(\de)\geq \delta \mathbf{1}_m$ for all $t\in[0, T]$.
\end{assumption}

A small investor, whose actions cannot affect the asset prices, will decide at every time
$t\in[0, T]$ the amount $\pi_{j,t}$ of his wealth to invest in the $j$-th risky asset, $j=1, \ldots, m$. The vector process $\pi:=(\pi_1, \ldots, \pi_m)^{\top}$ is called a portfolio of the investor. Then the investor's self-financing wealth process $X$ corresponding to a portfolio $\pi$ is the unique strong solution of the SDE:
\begin{align}
\label{wealth}
\begin{cases}
\dd X_t=[r_tX_{t-}+\pi_t^{\top}\mu_t]\dt+\pi_t^{\top}\sigma_t\dw_t+\int_{\cE}\pi_t^{\top}F_t(e) \widetilde N(\dt,\de), \\
X_0=x.
\end{cases}
\end{align}

The admissible portfolio set is defined as
\begin{align*}
\mathcal U=\Big\{\pi\in L^2_{\mathbb F}(0, T;\mathbb R^m)\;\Big|\; \pi_t\in\Pi~\dpt \Big\},
\end{align*}
where $\Pi\in\R^m$ is a given closed convex cone.
For instance, $\Pi=\R^{m}$ means there is no trading constraint; while $\Pi=\R_{+}^{m}$ means shorting is not allowed in the market.
For any $\pi\in \mathcal{U}$, the SDE \eqref{wealth} has a unique strong solution.
Different from the previous sections, in this section we request the constraint set $\Pi$ to be convex in order to apply the dual approach below.

\par
For a given expectation level $z\geq xe^{\int_0^Tr_s\ds}$, the investor's mean-variance problem is to
\begin{align}
\mathrm{Minimize}&\quad \mathrm{Var}(X_T)\equiv\E[X_T^{2}-z^2]%
, \nn\\
\mathrm{ s.t.} &\quad
\begin{cases}
\E[X_T]=z, \\
\pi\in \mathcal{U}.
\end{cases}
\label{optm}%
\end{align}
\begin{remark}
Lim \cite{Lim} studied a mean-variance problem with jumps without portfolio constraints, i.e. $\Pi=\R^m$. In his model, all the coefficients in \eqref{wealth} are assumed to be predictable with respect to the Brownian motion filtration, so no jump term has entered into his SRE, which is exactly the same one as in the model without jumps.
\end{remark}

We shall say that the mean-variance problem \eqref{optm} is feasible for a given level $z\geq xe^{\int_0^Tr_s\ds}$ if there is a portfolio $\pi\in \mathcal{U}$ which satisfies the target constraint $\E[X_T]=z$. An optimal portfolio to \eqref{optm} is called an efficient portfolio corresponding to $z$ and the corresponding $(\sqrt{\mathrm{Var}(X_T)},z)$ is called an efficient point. The set of all efficient points, with $z\geq xe^{\int_0^Tr_s\ds}$, is called the efficient frontier.

Define the dual cone of $\Pi$ as
\[\dualpi:=\Big\{y\in\R^{m}\;\Big|\; x^{\top}y\leq 0 \mbox{ for all $x\in\Pi$}\Big\}.\]
The following result gives an equivalent condition for the feasibility of \eqref{optm}. The proof is exactly the same as \cite[Theorem 5.3]{HSX}, so we omit it.
\begin{theorem}[Feasibility]\label{theoremfeasible}
Under assumption \ref{assu4}, the mean-variance problem \eqref{optm} is feasible for any $z\geq x e^{\int_0^Tr_t\dt}$ if and only if
\begin{equation}
\label{feasible}
\int_0^T \mathbb{P}(\mu_t\notin\dualpi) \dt>0.
\end{equation}
\end{theorem}
For the rest of this section, we will always assume \eqref{feasible} holds.

The way to solve \eqref{optm} is rather clear nowadays. To deal with the constraint $\E[X_T]=z$, we introduce a Lagrange
multiplier $-2\lambda\in\R$ and obtain the following \emph{relaxed}
optimization problem:
\begin{align}\label{optmun}
\inf &\quad J(x,\pi, \lambda; z)={\mathbb{E}}[(X_T-\lambda)^{2}]-(\lambda-z)^{2}, \\
\mathrm{s.t.} &\quad \pi\in \mathcal{U}.\nn%
\end{align}
Denote its optimal value as $$V(x,\lambda ; z)=\inf_{\pi\in\mathcal{U}}J(x,\pi,\lambda ; z).$$
According to the Lagrange duality theorem (see Luenberger \cite{Lu})
\begin{align}\label{duality}
\inf_{\pi\in\mathcal{U}, \E[X_T]=z}\mathrm{Var}(X_T%
)=\sup_{\lambda\in\R}V(x,\lambda ; z).
\end{align}
So we can solve the problem \eqref{optmun} by a two-step procedure: Firstly determine $V(x,\lambda ; z)$ for every $\lambda$, and then try to find a $\lambda^{*}$ to maximize $\lambda\mapsto V(x,\lambda;z)$.

\par

The relaxed problem \eqref{optmun} is a special stochastic LQ problem \eqref{LQ} studied in Section \ref{section:LQ}, where
\begin{align}\label{coeffMV}
A=r, \ B=\mu, \ C=0, \ D=\sigma^{\top}, \ E=0, \ Q=0, \ R=0, \ S=0, \ G=1.
\end{align}
The associated BSDEJ \eqref{P1} becomes
\begin{align}\label{Pmv}
\begin{cases}
\dd P_{1,t}=-\Big[2rP_{1,t-}+H_1^*(t,P_{1,t-},P_{2,t-},\Lambda_{1,t},\Gamma_{1,t},\Gamma_{2,t})\Big]\dt+\Lambda_{1,t}^{\top}\dw+\int_{\cE}\Gamma_{1,t}(e) \widetilde N(\dt,\de),\\
\dd P_{2,t}=-\Big[2rP_{2,t-}+H_2^*(t,P_{1,t-},P_{2,t-},\Lambda_{2,t},\Gamma_{1,t},\Gamma_{2,t})\Big]\dt+\Lambda_{2,t}^{\top}\dw+\int_{\cE}\Gamma_{2,t}(e) \widetilde N(\dt,\de),\\
P_{1,T}=1, \ P_{2,T}=1,
\end{cases}
\end{align}
where $H_1^*, \ H_2^*$, $\hat v_1$, $\hat v_2$ are defined as in \eqref{Hstar1}, \eqref{Hstar2} and \eqref{hatv} with coefficients given in \eqref{coeffMV}:
\begin{align*}
H_1(\omega,t,v,P_{1},P_{2},\Lambda_{1},\Gamma_{1},\Gamma_2)&=P_{1}v^{\top}\sigma\sigma^{\top}v+2(P_{1}\mu+\sigma\Lambda_{1})^{\top}v\\
&\qquad+\int_{\cE}\Big[(P_{1}+\Gamma_{1})\Big(((1+F^{\top}v)^+)^2-1\Big)-2P_{1}F^{\top}v\\
&\qquad\qquad+(P_{2}+\Gamma_{2})((1+F^{\top}v)^-)^2\Big]\nu(\de),\\
H_2(\omega,t,v,P_{1},P_{2},\Lambda_{2},\Gamma_{1},\Gamma_2)&=P_{2}v^{\top}\sigma\sigma^{\top}v-2(P_{2}\mu+\sigma\Lambda_{2})^{\top}v\\
&\qquad+\int_{\cE}\Big[(P_{2}+\Gamma_{2})\Big(((-1+F^{\top}v)^-)^2-1\Big)+2P_{2}F^{\top}v\\
&\qquad\qquad+(P_{1}+\Gamma_{1})((-1+F^{\top}v)^+)^2\Big]\nu(\de).
\end{align*}

Clearly, Theorems \ref{Th:singular} and \ref{uniqueness} can be applied to the BSDEJ \eqref{Pmv} to ensure that it admits a unique uniformly positive solution $(P_1,\Lambda_1,\Gamma_1, P_2,\Lambda_2,\Gamma_2)$.
Accordingly, Theorem \ref{Th:verif} leads to the following solution to the relaxed problem \eqref{optmun}.
\begin{theorem}
Let $(P_1,\Lambda_{1},\Gamma_1,P_2,\Lambda_{2},\Gamma_2)$ be the unique uniformly positive solution to \eqref{Pmv}. Then the state feedback control given by
\begin{align}\label{pistar}
\pi^*(t,X)&=\hat v_1(\omega,t,P_{1},P_{2},\Lambda_{1},\Gamma_{1},\Gamma_2)\Big(X_{t-}-\lambda e^{-\int_t^Tr_s\ds}\Big)^+\nn\\
&\qquad+\hat v_2(\omega,t,P_{1},P_{2},\Lambda_{2},\Gamma_{1},\Gamma_2)\Big(X_{t-}-\lambda e^{-\int_t^Tr_s\ds}\Big)^-,
\end{align}
is optimal for the LQ problem \eqref{LQ}. Moreover, the optimal value is
\begin{align*}
V(x,\lambda;z)=P_{1,0}\Big[\Big(x-\lambda e^{-\int_0^Tr_s\ds}\Big)^+\Big]^2+P_{2,0}\Big[\Big(x-\lambda e^{-\int_0^Tr_s\ds}\Big)^-\Big]^2-(\lambda-z)^2.
\end{align*}
\end{theorem}
This resolves the first step problem.
To solve the second step problem, i.e., to maximize $\lambda\mapsto V(x,\lambda;z)$, the following result is critical.
\begin{lemma}
Assume Assumption \ref{assu4} and conditioin \eqref{feasible} hold. Then
\begin{align}\label{strict}
P_{1,0}e^{-2\int_0^Tr_s\ds}-1\leq 0, \ P_{2,0}e^{-2\int_0^Tr_s\ds}-1< 0.
\end{align}
\end{lemma}
\pf
Applying It\^{o}'s formula to $P_{2,t}e^{-2\int_t^Tr_s\ds}$ on $[0,T]$, we have
\begin{align}\label{P2strict}
1-P_{2,0}e^{-2\int_0^Tr_s\ds}=
-\E\int_0^TH_2^*(t,P_{2,t-},\Lambda_{2,t},\Gamma_{2,t},P_{1,t-},\Gamma_{1,t})\dt.
\end{align}
Since $H_2^*(t,P_{2,t-},\Lambda_{2,t},\Gamma_{2,t},P_{1,t-},\Gamma_{1,t})\leq 0$ by its very definition, it follows $P_{2,0}e^{-2\int_0^Tr_s\ds}-1\leq 0$. Similarly, we can prove that $P_{1,0}e^{-2\int_0^Tr_s\ds}-1\leq 0$.

It remains to prove the strict inequality $P_{2,0}e^{-2\int_0^Tr_s\ds}-1< 0$. Suppose, on the contrary, $P_{2,0}e^{-2\int_0^Tr_s\ds}-1=0$. It then follows from \eqref{P2strict} that $H_2^*(t,P_{1,t-},P_{2,t-},\Lambda_{2,t},\Gamma_{1,t},\Gamma_{2,t})=0~\dpt$. Thus we deduce, from the uniqueness (Theorem \ref{uniqueness}) of solution to the BSDE \eqref{Pmv}, that $P_{2,t}=e^{2\int_t^Tr_s\ds}$, $\Lambda_{2,t}=0$ and $\Gamma_{2,t}=0$.

On the other hand,
\begin{align}\label{Fpositive}
((1-F^{\top}v)^+)^2 +2F^{\top}v-1\leq (1-F^{\top}v)^2 +2F^{\top}v-1=|F^{\top}v|^2.
\end{align}
Since $F\in L_{\mathcal{P}}^{\nu,\infty}(0,T;\R^m)$, there exists $c_1>0$ such that $|F(e)|\leq c_1$ for almost all $e\in\cE$. Hence
\begin{align}\label{Fnegative}
1-F^{\top}v\geq 0, \ \mbox{if} \ v\in\Pi \ \mbox{and } \ |v|\leq c_1^{-1}.
\end{align}
Combining \eqref{Fpositive} and \eqref{Fnegative}, we have
\begin{align*}
H_2^*(t,P_{1},P_{2},0,\Gamma_1 ,0)&=\inf_{v\in\Pi}\Big[P_{2}v^{\top}\sigma\sigma^{\top}v-2 P_{2}\mu^{\top}v+\int_{\cE}\Big[P_{2}\Big(((1-F^{\top}v)^+)^2 +2F^{\top}v-1\Big)\\
&\qquad\qquad+(P_{1}+\Gamma_{1})((1-F^{\top}v)^-)^2\Big]\nu(\de)\Big]\\
&\leq\inf_{v\in\Pi,|v|\leq c_1^{-1}}\Big[P_{2}v^{\top}\Big(\sigma\sigma^{\top}+\int_{\cE}FF^{\top}\nu(\de)\Big)v-2 P_{2}\mu^{\top}v\Big]\\
&\leq P_2\inf_{v\in\Pi,|v|\leq c_1^{-1}}\Big(c|v|^2-2\mu^{\top}v\Big).
\end{align*}
Note \eqref{feasible} implies that there exists $\mathcal{O}\subseteq \Omega\times[0,T]$ such that $\mu\notin\dualpi~\dpt$ on $\mathcal{O}$. Hence, there exists $v_0\in\Pi$ such that $\mu^{\top}v_0>0$ on $\mathcal{O}$. By choosing $v_1=\varepsilon v_0$ with $\varepsilon>0$ being sufficiently small so that $|v_1|\leq c_1^{-1}$, we get
\begin{align*}
H_2^*(t,P_{1},P_{2},0,\Gamma_1 ,0)
&\leq P_2\varepsilon\Big(c\varepsilon|v_0|^2-2\mu^{\top}v_0\Big)~\textrm{on $\mathcal{O}$.}
\end{align*}
The RHS is negative for sufficiently small $\varepsilon>0$ on $\mathcal{O}$, leading to a contraction. Therefore $P_{2,0}e^{-2\int_0^Tr_s\ds}-1< 0$.
\eof

To find a $\lambda^{*}$ to maximize $\lambda\mapsto V(x,\lambda;z)$, we do some tedious calculation (using \eqref{strict}) and obtain
$$\max_{\lambda}V(x,\lambda;z)= V(x,\lambda^*;z)=\frac{P_{2,0}}{1-P_{2,0}e^{-2\int_0^Tr_s\ds}}\Big(x-ze^{-\int_0^Tr_s\ds}\Big)^2,$$
where
$$\lambda^*=\frac{z-xP_{2,0}e^{-\int_0^Tr_s\ds}}{1-P_{2,0}e^{-2\int_0^Tr_s\ds}}.$$

The above analysis boils down to the following solution to the mean-variance problem \eqref{optm}.
\begin{theorem}
Let $(P_1,\Lambda_{1},\Gamma_1,P_2,\Lambda_{2},\Gamma_2)$ be the unique uniformly positive solution to \eqref{Pmv}. Then the state feedback portfolio given by
\begin{align}\label{efficient}
\pi^*(t,X)&=\hat v_1(\omega,t,P_{1},P_{2},\Lambda_{1},\Gamma_1,\Gamma_2)\Big(X_{t-}-\lambda^* e^{-\int_t^Tr_s\ds}\Big)^+\nn\\
&\qquad+\hat v_2(\omega,t,P_{1},P_{2},\Lambda_{1},\Gamma_2,\Gamma_2)\Big(X_{t-}-\lambda^* e^{-\int_t^Tr_s\ds}\Big)^-,
\end{align}
is optimal to the mean-variance problem \eqref{optm}. Moreover, the efficient frontier is determined by
\begin{align*}
\mathrm{Var}(X_T)=\frac{P_{2,0}e^{-2\int_0^Tr_s\ds}}{1-P_{2,0}e^{-2\int_0^Tr_s\ds}}\Big(\E[X_T]-xe^{\int_0^Tr_s\ds}\Big)^2,
\end{align*}
where $\E[X_T]\geq xe^{\int_0^Tr_s\ds}$.
\end{theorem}

\begin{remark}
In the constrained mean-variance model without jumps studied in Hu and Zhou \cite{HZ}, the efficient portfolio only takes the second term on the RHS of \eqref{efficient}, i.e. the optimal wealth $X_t$ will never exceed $\lambda^* e^{-\int_t^Tr_s\ds}$ on $[0,T]$, and it only depends on $(P_2,\Lambda_2)$ as $\hat v_2$ does. But in our cone-constrained MV problem with jumps, the optimal wealth $X_t$ will probably cross the bliss point $\lambda^* e^{-\int_t^Tr_s\ds}$.
\end{remark}

\begin{appendix}
\section{Heuristic derivation of the BSDEJ \eqref{P1}}\label{appnA}
By the Meyer-It\^o formula \cite[Theorem 70]{Protter}, we have
\begin{align*}
\dd X_t^+&=\mathbf{1}_{\{X_{t-}>0\}}\Big[\big(A_tX_{t-}+B_t^{\top}u_t-\int_{\cE}(E_t(e)X_{t-}+F_t(e)^{\top}u_t)\nu(\de)\big)\dt
+(C_tX_{t-}+D_tu_t)^{\top}\dw_t\Big]\\
&\qquad+\int_{\cE}\Big[(X_{t-}+E_t(e)X_{t-}+F_t(e)^{\top}u_t)^+-X_{t-}^+\Big]N(\dt,\de)+\frac{1}{2}\dd L_t,
\end{align*}
and
\begin{align*}
\dd X_t^-&=-\mathbf{1}_{\{X_{t-}\leq0\}}\Big[\big(A_tX_{t-}+B_t^{\top}u_t-\int_{\cE}(E_t(e)X_{t-}+F_t(e)^{\top}u_t)\nu(\de)\big)\dt
+(C_tX_{t-}+D_tu_t)^{\top}\dw_t\Big]\\
&\qquad+\int_{\cE}\Big[(X_{t-}+E_t(e)X_{t-}+F_t(e)^{\top}u_t)^--X_{t-}^-\Big]N(\dt,\de)+\frac{1}{2}\dd L_t,
\end{align*}
where $L$ is the local time of $X$ at $0$. Since $X_t^{\pm}\dd L_t=0$,
applying the It\^{o} formula yields
\begin{align*}
\dd\;(X_t^+)^2&=\mathbf{1}_{\{X_{t-}>0\}}\Big[2X_{t-}^+\big(A_tX_{t-}+B_t^{\top}u_t-\int_{\cE}(E_t(e)X_{t-}+F_t(e)^{\top}u_t)\nu(\de)\big)
+|C_tX_{t-}+D_tu_t|^2\Big]\dt\\
&\qquad+ \mathbf{1}_{\{X_{t-}>0\}}2X_{t-}^+(C_tX_{t-}+D_tu_t)^{\top}\dw_t\\
&\qquad+\int_{\cE}\Big[((X_{t-}+E_t(e)X_{t-}+F_t(e)^{\top}u_t)^+)^2-(X_{t-}^+)^2\Big]N(\dt,\de),
\end{align*}
and
\begin{align*}
\dd\;(X_t^-)^2&=\mathbf{1}_{\{X_{t-}\leq0\}}\Big[-2X_{t-}^-\big(A_tX_{t-}+B_t^{\top}u_t-\int_{\cE}(E_t(e)X_{t-}+F_t(e)^{\top}u_t)\nu(\de)\big)
+|C_tX_{t-}+D_tu_t|^2\Big]\dt\\
&\qquad- \mathbf{1}_{\{X_{t-}\leq0\}}2X_{t-}^-(C_tX_{t-}+D_tu_t)^{\top}\dw_t\\
&\qquad+\int_{\cE}\Big[((X_{t-}+E_t(e)X_{t-}+F_t(e)^{\top}u_t)^-)^2-(X_{t-}^-)^2\Big]N(\dt,\de).
\end{align*}
Assume that $P_1$ and $P_2$ are semimartingales of the following form:
\begin{align*}
\begin{cases}
\dd P_{1}=-f_1\dt+\Lambda_{1}^{\top}\dw+\int_{\cE}\Gamma_{1}(e) \widetilde N(\dt,\de),\\
P_{1,T}=G,
\end{cases}
\end{align*}
and
\begin{align*}
\begin{cases}
\dd P_{2}=-f_2\dt+\Lambda_{2}^{\top}\dw+\int_{\cE}\Gamma_{2}(e) \widetilde N(\dt,\de),\\
P_{2,T}=G.
\end{cases}
\end{align*}
Applying It\^{o}'s formula to $P_{1,t}(X_t^+)^2$,
\begin{align*}
\dd P_{1,t}(X_t^+)^2&=P_{1,t-}\mathbf{1}_{\{X_{t-}>0\}}\Big[2X_{t-}^+(A_tX_{t-}+B_t^{\top}u_t)+|C_tX_{t-}+D_tu_t|^2\Big]\dt\\
&\qquad+2X_{t-}^+(C_tX_{t-}+D_tu_t)^{\top}\Lambda_{1,t}\dt\\
&\qquad-2P_{1,t-}X_{t-}^+\mathbf{1}_{\{X_{t-}>0\}}\int_{\cE}(E_t(e)X_{t-}+F_t(e)^{\top}u_t)\nu(\de)\dt\\
&\qquad+\int_{\cE}(P_{1,t-}+\Gamma_{1,t}(e))\Big(((X_{t-}+E_t(e)X_{t-}+F_t(e)^{\top}u_t)^+)^2-(X_{t-}^+)^2\Big)\nu(\de)\dt\\
&\qquad-(X_{t-}^+)^2f_1\dt+\Big[(X_{t-}^+)2\Lambda_{1,t}+2P_{1,t-}X_{t-}^+(C_tX_{t-}+D_tu_t)\Big]^{\top}\dw\\
&\qquad+\int_{\cE} \Big[(P_{1,t-}+\Gamma_{1}(e))[((X_{t-}+E_t(e)X_{t-}+F_t(e)^{\top}u_t)^+)^2-(X_{t-}^+)^2]\\
&\qquad\qquad\qquad+(X_{t-}^+)^2\Gamma_{1}(e) \Big] \widetilde N(\dt,\de),
\end{align*}
and
applying It\^{o}'s formula to $P_{2,t}(X_t^-)^2$,
\begin{align*}
\dd P_{2,t}(X_t^-)^2&=P_{2,t-}\mathbf{1}_{\{X_{t-}\leq0\}}\Big[-2X_{t-}^-(A_tX_{t-}+B_t^{\top}u_t)+|C_tX_{t-}+D_tu_t|^2\Big]\dt\\
&\qquad-2X_{t-}^-(C_tX_{t-}+D_tu_t)^{\top}\Lambda_{2,t}\dt\\
&\qquad+2P_{2,t-}X_{t-}^-\mathbf{1}_{\{X_{t-}\leq0\}}\int_{\cE}(E_t(e)X_{t-}+F_t(e)^{\top}u_t)\nu(\de)\dt\\
&\qquad+\int_{\cE}(P_{2,t-}+\Gamma_{2,t}(e))\Big(((X_{t-}+E_t(e)X_{t-}+F_t(e)^{\top}u_t)^-)^2-(X_{t-}^-)^2\Big)\nu(\de)\dt\\
&\qquad-(X_{t-}^-)^2f_2\dt+\Big[(X_{t-}^-)^2\Lambda_{2,t}-2P_{2,t-}X_{t-}^-(C_tX_{t-}+D_tu_t)\Big]^{\top}\dw\\
&\qquad+\int_{\cE} \Big[(P_{2,t-}+\Gamma_{2}(e))[((X_{t-}+E_t(e)X_{t-}+F_t(e)^{\top}u_t)^-)^2-(X_{t-}^-)^2]\\
&\qquad\qquad\qquad+(X_{t-}^-)^2\Gamma_{2}(e) \Big] \widetilde N(\dt,\de).
\end{align*}
Then
\begin{align*}
J(x; u)&=\mathbb{E}\Big[\int_0^T\Big(Q_tX_t^2+u_t^{\top}R_tu_t+2X_tS_t^{\top}u_t\Big)\dt+GX_T^2\Big]\\
&=P_{1,0}(x^+)^2+P_{2,0}(x^-)^2+\mathbb{E}\int_0^T\Big(Q_tX_t^2+u_t^{\top}R_tu_t+2X_tS_t^{\top}u_t\Big)\dt\\
&\qquad+P_{1,t-}\mathbf{1}_{\{X_{t-}>0\}}\Big[2X_{t-}^+(A_tX_{t-}+B_t^{\top}u_t)+|C_tX_{t-}+D_tu_t|^2\Big]\dt\\
&\qquad+2X_{t-}^+(C_tX_{t-}+D_tu_t)^{\top}\Lambda_{1,t}\dt\\
&\qquad-2P_{1,t-}X_{t-}^+\mathbf{1}_{\{X_{t-}>0\}}\int_{\cE}(E_t(e)X_{t-}+F_t(e)^{\top}u_t)\nu(\de)\dt\\
&\qquad+\int_{\cE}(P_{1,t-}+\Gamma_{1,t}(e))\Big(((X_{t-}+E_t(e)X_{t-}+F_t(e)^{\top}u_t)^+)^2-(X_{t-}^+)^2\Big)\nu(\de)\dt\\
&\qquad-(X_{t-}^+)^2f_1\dt-(X_{t-}^-)^2f_2\dt
\end{align*}
\begin{align*}
&\qquad +P_{2,t-}\mathbf{1}_{\{X_{t-}\leq0\}}\Big[-2X_{t-}^-(A_tX_{t-}+B_t^{\top}u_t)+|C_tX_{t-}+D_tu_t|^2\Big]\dt\\
&\qquad-2X_{t-}^-(C_tX_{t-}+D_tu_t)^{\top}\Lambda_{2,t}\dt\\
&\qquad+2P_{2,t-}X_{t-}^-\mathbf{1}_{\{X_{t-}\leq0\}}\int_{\cE}(E_t(e)X_{t-}+F_t(e)^{\top}u_t)\nu(\de)\dt\\
&\qquad+\int_{\cE}(P_{2,t-}+\Gamma_{2,t}(e))\Big(((X_{t-}+E_t(e)X_{t-}+F_t(e)^{\top}u_t)^-)^2-(X_{t-}^-)^2\Big)\nu(\de)\dt.
\end{align*}
Denote $\phi(X_{t-},u_t)$ be the integrand on the RHS of the above equation.
\begin{itemize}
\item If $X_{t-}>0$, then
\begin{align*}
\phi(X_{t-},u_t)&=\Big[Q+v^{\top}(R+P_{1,t-}D^{\top}D)v+2S^{\top}v+P_{1,t-}(2A+C^{\top}C)+2C^{\top}\Lambda_{1,t}\\
&\qquad+2(P_{1,t-}(B+D^{\top}C)+D^{\top}\Lambda_1)^{\top}v\\
&\qquad-2P_{1,t-}\int_{\cE}(E+F^{\top}v)\nu(\de)
+\int_{\cE}(P_{1,t-}+\Gamma_{1,t-})\Big(((1+E+F^{\top}v)^+)^2-1\Big)\nu(\de)\\
&\qquad-f_1+\int_{\cE}(P_{2,t-}+\Gamma_{2,t-})((1+E+F^{\top}v)^-)^2\nu(\de)\Big]X_{t-}^2,
\end{align*}
where $v_t=\frac{u_t}{|X_{t-}|}$.

\item
If $X_{t-}<0$, then
\begin{align*}
\phi(X_{t-},u_t)&=\Big[Q+v^{\top}(R+P_{2,t-}D^{\top}D)v-2S^{\top}v+P_{2,t-}(2A+C^{\top}C)+2C^{\top}\Lambda_{2,t}\\
&\qquad-2(P_{1,t-}(B+D^{\top}C)+D^{\top}\Lambda_2)^{\top}v\\
&\qquad-2P_{2,t-}\int_{\cE}(E-F^{\top}v)\nu(\de)
+\int_{\cE}(P_{2,t-}+\Gamma_{2,t-})\Big(((-1-E+F^{\top}v)^-)^2-1\Big)\nu(\de)\\
&\qquad-f_2+\int_{\cE}(P_{1,t-}+\Gamma_{1,t-})((-1-E+F^{\top}v)^+)^2\nu(\de)\Big]X_{t-}^2,
\end{align*}
where $v_t=\frac{u_t}{|X_{t-}|}$.

\item
If $X_{t-}=0$, then $\phi(0, 0)=0$ and
\begin{align*}
\phi(X_{t-},u_t)&=u_t^{\top}Ru_t+P_{2,t-}|Du_t|^2
+\int_{\cE}(P_{1,t-}+\Gamma_{1,t})((F^{\top}u_t)^+)^2\nu(\de)\\
&\qquad+\int_{\cE}(P_{2,t}+\Gamma_{2,t-})((F^{\top}u_t)^-)^2\nu(\de)\geq 0.
\end{align*}

\end{itemize}

In order to ensure $\min_{u\in\Pi}\phi(X_t,u)=0$, it is evident to take $f_1$ and $f_2$ in the form of \eqref{Hstar1} and \eqref{Hstar2}.

\section{Proof of Lemma \ref{thlemma1}}\label{appnB}

For any positive integers $k< l$, set
\[
P_{i,t}^{k,l}:=P_{i,t}^{k}-P_{i,t}^{l}\geq0, \ \Lambda_{i,t}^{k,l}:=\Lambda_{i,t}^{k}-\Lambda_{i,t}^{l}, \ \Gamma_{i,t}^{k,l}:=\Gamma_{i,t}^{k}-\Gamma_{i,t}^{l}, \ i=1,2.
\]
Let $\kappa>0$ be a constant to be specified later, and write $$\Psi(x)=\frac{1}{\kappa}\big(e^{\kappa x}-\kappa x-1\big).$$
Notice $\Psi'(x)\geq 0$ for $x\geq 0$.
Applying It\^{o}'s formula to $\Psi(P_{1,t}^{k,l})$, we get
\begin{align*}
&\quad\Psi(P_{1,0}^{k,l})+\frac{1}{2}\E\int_0^T\Psi''(P_{1,t}^{k,l})|\Lambda_1^{k,l}|^2\dt\\
&\qquad\qquad+\E\int_0^T\int_{\cE}\Big[\Psi(P_{1,t-}^{k,l}+\Gamma_{1,t}^{k,l})-\Psi(P_{1,t-}^{k,l})
-\Psi'(P_{1,t-}^{k,l})\Gamma_1^{k,l}\Big]\nu(\de)\dt\\
&=\Psi(0)+\E\int_0^T\Psi'(P_{1,t}^{k,l})\Big[(2A+C^{\top}C)P^{k,l}_{1,t-}+2C^{\top}\Lambda^{k,l}_{1,t}\\
&\qquad+H_1^{k}(t,P_{1,t-}^k,P_{2,t-}^k,\Lambda_{1,t}^k,\Gamma_{1,t}^k,\Gamma_{2,t}^k)
-H_1^{l}(t,P_{1,t-}^l,P_{2,t-}^l,\Lambda_{1,t}^l,\Gamma_{1,t}^l,\Gamma_{2,t}^l)\Big]\dt.
\end{align*}
Using the following fact:
\begin{align*}
&\Psi(0)=0, \ P_{1,t}^{k,l}\geq 0, \ \Psi'(P_{1,t}^{k,l})=e^{\kappa P_{1,t}^{k,l}}-1\geq 0,\ H_1^{l}\geq H_1^{*}, \\
&H_1^{k}\leq \int_{\cE}\Big[(P_{1}+\Gamma_{1})\Big(((1+E)^+)^2-1\Big)-2P_{1}E+(P_{2}+\Gamma_{2})((1+E)^-)^2\Big]\nu(\de)\leq c,
\end{align*}
where $c$ is independent of $k$ and $l$,
we obtain
\begin{align*}
&\quad\Psi(P_{1,0}^{k,l})+\frac{1}{2}\E\int_0^T\Psi''(P_{1,t}^{k,l})|\Lambda_1^{k,l}|^2\dt\\
&\qquad\qquad+\E\int_0^T\int_{\cE}\Big[\Psi(P_{1,t-}^{k,l}+\Gamma_{1,t}^{k,l})-\Psi(P_{1,t-}^{k,l})
-\Psi'(P_{1,t-}^{k,l})\Gamma_1^{k,l}\Big]\nu(\de)\dt\\
&\leq\E\int_0^T\Psi'(P_{1,t}^{k,l})\Big[c+2C^{\top}\Lambda_{1,t}^{k,l}
-H_1^{*}(t,P_{1,t-}^l,P_{2,t-}^l,\Lambda_{1,t}^l,\Gamma_{1,t}^l,\Gamma_{2,t}^l)\Big]\dt.
\end{align*}
Keeping in mind $ P^l_i$, $\Gamma^l_i$, $ i=1,2$, are uniformly bounded, we have the following estimates:
\begin{align*}
-H_1^{*}(t,P_{1,t-}^l,P_{2,t-}^l,\Lambda_{1,t}^l,\Gamma_{1,t}^l,\Gamma_{2,t}^l)&\leq -\inf_{v\in\R^m}H_1(t,0,P_{1,t-}^l,P_{2,t-}^l,\Lambda_{1,t}^l,\Gamma_{1,t}^l,\Gamma_{2,t}^l)\\
&\leq c+c|\Lambda_1^l|^2\\
&\leq c+3c(|\Lambda_1^{k,l}|^2+|\Lambda_1^{k}-\Lambda_1|^2+|\Lambda_1|^2),
\end{align*}
where $c>0$ is a constant independent of $l$ and $k$.
The above estimates lead to
\begin{align}
&\quad\Psi(P_{1,0}^{k,l})+\E\int_0^T\Big(\frac{1}{2}\Psi''(P_{1,t}^{k,l})-3c\Psi'(P_{1,t}^{k,l})\Big) |\Lambda_1^{k,l}|^2\dt\nn\\
&\qquad\qquad+\E\int_0^T\int_{\cE}\Big[\Psi(P_{1,t-}^{k,l}+\Gamma_{1,t}^{k,l})-\Psi(P_{1,t-}^{k,l})
-\Psi'(P_{1,t-}^{k,l})\Gamma_1^{k,l}\Big]\nu(\de)\dt\nn\\
&\leq\E\int_0^T\Psi'(P_{1,t}^{k,l})\Big[2c+2C^{\top}\Lambda_{1,t}^{k,l}+3c|\Lambda_1^{k}-\Lambda_1|^2+3c|\Lambda_1|^2)
\Big]\dt.\label{estimates1}
\end{align}

Take $\kappa=12c$. Then $\frac{1}{2}\Psi''(x)- 3c\Psi'(x)=3c(e^{\kappa x}+1)=3c\Psi'(x)+6c\geq 6c$ for $x\geq 0$.
So by the dominated convergence theorem, the sequence
$$\sqrt{\frac{1}{2}\Psi''(P_{1,t}^{k,l})-3c\Psi'(P_{1,t}^{k,l})}$$
converges strongly to
$$\sqrt{3c\Psi' (P_{1,t}^{k}-P_{1,t})+6c},$$
as $l\to \infty$, and they are uniformly bounded. Therefore,
$$\sqrt{3c\Psi'(P_{1,t}^{k,l})+6c}\;\Lambda_1^{k,l}$$
converges weakly to
$$\sqrt{3c\Psi' (P_{1,t}^{k}-P_{1,t})+6c}\; (\Lambda_1^k-\Lambda_1).$$
By the mean value theorem and the uniformly boundedness of $P_{1,t-}^{k,l}$ and $\Gamma_{1,t-}^{k,l}$, we obtain
\begin{align}\label{AGamma}
\Psi(P_{1,t-}^{k,l}+\Gamma_{1,t}^{k,l})-\Psi(P_{1,t-}^{k,l})
-\Psi'(P_{1,t-}^{k,l})\Gamma_1^{k,l}&=\frac{1}{\kappa}e^{\kappa P_{1,t-}^{k,l}}
\Big[e^{\kappa\Gamma_{1,t}^{k,l}}-\kappa\Gamma_{1,t}^{k,l}-1\Big] \geq \ep|\Gamma_{1,t}^{k,l}|^2,
\end{align}
for some constant $\ep>0$ independent of $l$ and $k$.
We then get from \eqref{estimates1} that
\begin{align*}
&\quad\E\int_0^T\big(3c\Psi' (P_{1}^{k}-P_{1})+6c\big)|\Lambda_1^{k}-\Lambda_1|^2\dt\\
&\leq\varliminf_{l\to \infty}\E\int_0^T\big(3c\Psi' (P_{1}^{k,l})+6c\big)|\Lambda_1^{k,l}|^2\dt\\
&\leq\E\int_0^T\Psi'(P_{1}^{k}-P_{1})\Big[2c+2C^{\top}(\Lambda_{1}^{k}-\Lambda_{1})+3c|\Lambda_1^{k}-\Lambda_1|^2+3c|\Lambda_1|^2)\Big]\ds.
\end{align*}
Canceling the common terms, it yields
\begin{align*}
\E\int_0^T6c|\Lambda_1^{k}-\Lambda_1|^2\ds \leq\E\int_0^T\Psi'(P_{1}^{k}-P_{1})\Big[2c+2C^{\top}(\Lambda_{1}^{k}-\Lambda_{1})+3c|\Lambda_1|^2)
\Big]\ds.
\end{align*}
By passing to the limit $k\to \infty$, applying dominated convergence theorem and noticing $\Psi'(0)=0$, we have
\begin{align*}
\lim_{k\to \infty}\E\int_0^T|\Lambda_1^k-\Lambda_1|^2\dt=0.
\end{align*}

Using \eqref{AGamma}, we can similarly get
\begin{align*}
\lim_{k\to \infty}\E\int_0^T\int_{\cE}|\Gamma_1^k-\Gamma_1|^2\nu(\de)\dt=0.
\end{align*}
Because $\Gamma_1^k$, $k=1,2,\cdots,$ are uniformly bounded in $\linnu(0, T;\R) $, we conclude that $\Gamma_1 \in \linnu(0, T;\R) $.

Along appropriate subsequence (which is still denoted by $(\Lambda_{1}^k, \Gamma_{1}^k)$) we may obtain $\dpt$ convergence of
\begin{align*}
\int_t^T(\Lambda_{1}^k)^{\top}\dw+\int_t^T\int_{\cE}\Gamma_1^k(e) \widetilde N(\ds,\de)\to \int_t^T \Lambda_{1}^{\top}\dw+\int_t^T\int_{\cE}\Gamma_1(e) \widetilde N(\ds,dy),
\end{align*}
and
\begin{align*}
&\lim_{k\to \infty}\int_t^T\Big[(2A+C^{\top}C)P_{1,t-}^k+2C^{\top}\Lambda_{1,t}^k+Q\Big]\ds
=\int_t^T \Big[(2A+C^{\top}C)P_{1,t-}+2C^{\top}\Lambda_{1,t}+Q\Big]\ds.
\end{align*}
We now turn to prove
\begin{align}\label{AFG}
&\lim_{k\to \infty}\int_t^TH_1^{k}(s,P_{1,s-}^k,P_{2,s-}^k,\Lambda_{1,s}^k,\Gamma_{1,s}^k,\Gamma_{2,s}^k)\ds
=\int_t^T H_1^{*}(s,P_{1,s-},P_{2,s-},\Lambda_{1,s},\Gamma_{1,s},\Gamma_{2,s})\ds.
\end{align}
We have
\begin{align*}
&\quad|H_1^{k}(s,P_{1,s-}^k,P_{2,s-}^k,\Lambda_{1,s}^k,\Gamma_{1,s}^k,\Gamma_{2,s}^k)-H_1^{*}(s,P_{1,s-},P_{2,s-},\Lambda_{1,s},\Gamma_{1,s},\Gamma_{2,s})|\\
&\leq |H_1^{k}(s,P_{1,s-}^k,P_{2,s-}^k,\Lambda_{1,s}^k,\Gamma_{1,s}^k,\Gamma_{2,s}^k)-H_1^{*}(s,P_{1,s-}^k,P_{2,s-}^k,\Lambda_{1,s}^k,\Gamma_{1,s}^k,\Gamma_{2,s}^k)|\\
&\qquad +|H_1^{*}(s,P_{1,s-}^k,P_{2,s-}^k,\Lambda_{1,s}^k,\Gamma_{1,s}^k,\Gamma_{2,s}^k)-H_1^{*}(s,P_{1,s-},P_{2,s-},\Lambda_{1,s},\Gamma_{1,s},\Gamma_{2,s})|.
\end{align*}
Recall that $\Lambda_{1,s}^k\to \Lambda_{1,s}~\dpt$, so there exists $k_{1}(\omega,s)$ such that $|\Lambda_{1,s}^k|\leq 1 +|\Lambda_{1,s}|$ for $k\geq k_{1}$.
Notice that $H_1(s,0,P_{1,s-}^k,\Lambda_{1,s}^k,\Gamma_{1,s}^k,P_{2,s-}^k+\Gamma_{2,s}^k)$ is upper bounded by some $c>0$, and
\begin{align*}
H_1(s,v,P_{1,s-}^k,P_{2,s-}^k,\Lambda_{1,s}^k,\Gamma_{1,s}^k,\Gamma_{2,s}^k)
&\geq \delta|v|^2-2c_1|v|(1+|\Lambda_1^k|)-c_1\\
&\geq \delta|v|^2-2c_1|v|(2+|\Lambda_1|)-c_1\geq c,
\end{align*}
if $|v|>c_2(1 +|\Lambda_{1,s}|)$ with $c_2>0$ being sufficiently large. Hence, for $k\geq\max\{ c_2(1+|\Lambda_{1,s}|), k_{1}\}$, we have
\begin{align*}
H_1^{*}(s,P_{1,s-}^k,P_{2,s-}^k,\Lambda_{1,s}^k,\Gamma_{1,s}^k,\Gamma_{2,s}^k)
&=\inf_{\substack{v\in\Pi\\|v|\leq c_2(1 +|\Lambda_{1}|)}}H_1(s,v,P_{1,s-}^k,P_{2,s-}^k,\Lambda_{1,s}^k,\Gamma_{1,s}^k,\Gamma_{2,s}^k)\\
&\geq \inf_{\substack{v\in\Pi\\|v|\leq k}}H_1(s,v,P_{1,s-}^k,P_{2,s-}^k,\Lambda_{1,s}^k,\Gamma_{1,s}^k,\Gamma_{2,s}^k)\\
&=H_1^{k}(s,P_{1,s-}^k,P_{2,s-}^k,\Lambda_{1,s}^k,\Gamma_{1,s}^k,\Gamma_{2,s}^k).
\end{align*}
We also have the reverse inequality $H_1^*\leq H_1^{k}$ by definition.
Therefore, when $k\geq\max\{ c_2(1+|\Lambda_{1,s}|), k_{1}\}$, $H_1^{k}$ and $H_1^{*}$ coincide at $(s,P_{1,s-}^k,P_{2,s-}^k,\Lambda_{1,s}^k,\Gamma_{1,s}^k,\Gamma_{2,s}^k)$.

Notice that
\begin{align*}
&\quad|H_1^{*}(s,P_{1,s-}^k,P_{2,s-}^k,\Lambda_{1,s}^k,\Gamma_{1,s}^k,\Gamma_{2,s}^k)
-H_1^{*}(s,P_{1,s-},P_{2,s-},\Lambda_{1,s},\Gamma_{1,s},\Gamma_{2,s}))|\\
&\leq \sup_{\substack{v\in\Pi\\|v|\leq c_2(1 +|\Lambda_{1}|)}}|H_1(s,v,P_{1,s-}^k,P_{2,s-}^k,\Lambda_{1,s}^k,\Gamma_{1,s}^k,\Gamma_{2,s}^k)
-H_1(s,v,P_{1,s-},P_{2,s-},\Lambda_{1,s},\Gamma_{1,s},\Gamma_{2,s}))|,
\end{align*}
one can easily see, from the definition of $H_1$, that as long as
$$|P_1^k-P_1|+|\Lambda_1^k-\Lambda_1|+\int_{\cE}|\Gamma_1^k-\Gamma_1|\nu(\de)+|P_2^k-P_2|+\int_{\cE}|\Gamma_2^k-\Gamma_2|\nu(\de)\to 0,$$
we have
$$\lim_{k\to \infty}|H_1^{*}(s,P_{1,s-}^k,P_{2,s-}^k,\Lambda_{1,s}^k,\Gamma_{1,s}^k,\Gamma_{2,s}^k)-H_1^{*}(s,P_{1,s-},P_{2,s-},\Lambda_{1,s},\Gamma_{1,s},\Gamma_{2,s}))|=0.$$


Since
\begin{align*}
|H_1^{k}(s,P_{1,s-}^k,\Lambda_{1,s}^k,\Gamma_{1,s}^k,P_{2,s-}^k,\Gamma_{2,s}^k)|\leq c(1+|\Lambda_{1,s}^k|^2),
\end{align*}
the dominated convergence theorem leads to \eqref{AFG}.

Now it is standard to show that $$\lim_{k\to\infty}\E\Big[\sup_{t\in[0,T]}|P^k_{1,t}-P_{1,t}|\Big]=0;$$ please refer to Antonelli and Mancini \cite[Theorem 1]{AM} for details.

\end{appendix}

%

\end{document}